\documentclass[12pt]{article}
\usepackage{amsmath}
\usepackage{amssymb}
\usepackage[ctr1,eng]{myarticl}

\newcounter{mysubs}[section]
\renewcommand{\themysubs}{\thesection.\arabic{mysubs}}
\newcommand{\mysubs}{
\refstepcounter{mysubs}
\addvspace{\baselineskip}\noindent \textbf{\themysubs} }

\newcommand{\cB }{\mathcal{B}}

\newcommand{\cG }{\mathcal{G}}
\newcommand{\cI }{\mathcal{I}}

\newcommand{\copr }{\varDelta }
\newcommand{\End }{\mathrm{End}}
\newcommand{\hatt }{\hat{t}}
\newcommand{\hatu }{\hat{u}}

\newcommand{\hatw }{\hat{w}}
\newcommand{\hatz }{\hat{z}}

\newcommand{\lact }{.}
\newcommand{\lcoa }{\delta }

\newcommand{\Ndbasis }{\boldsymbol{\mathrm{e}}}
\newcommand{\ndN }{\mathbb{N}}

\newcommand{\ndZ }{\mathbb{Z}}
\newcommand{\ot }{\otimes }
\newcommand{\overseteqnum }[2]{\overset{\makebox[0pt]{\tiny (\ref{#1})}}{#2}}
\newcommand{\paar }[2]{\langle #1,#2\rangle }
\newcommand{\paarb }[2]{\bigg\langle #1,#2\bigg\rangle }
\newcommand{\psubq }{\mbox{$\subset \!\!\!|$}}
\newcommand{\qfact }[2]{(#1)^!_{#2}}
\newcommand{\qnum }[2]{(#1)_{#2}}
\newcommand{\restr }[1]{\!\upharpoonright _{#1}}

\newcommand{\subq }{\mbox{$\subseteq \!\!\!|$}}

\newcommand{\YD }{Yetter--Drinfel'd }
\newcommand{\YDcat }{\mathcal{YD}}

\title{Finite dimensional rank 2 Nichols algebras of diagonal type II:
Classification}
\author{I.~Heckenberger\thanks{email:
Istvan.Heckenberger@math.uni-leipzig.de}}

\begin{document}
\maketitle

\begin{abstract}
The method of subquotients is developed and used to determine
all finite dimensional rank 2 Nichols algebras of diagonal type
over an arbitrary field of characteristic zero.


\end{abstract}

\section{Introduction}

This paper is the continuation of \cite{p-Heck04a}.
Let $k$ be a field of characteristic zero,
$G$ an abelian group, $V\in ^{kG}_{kG}\YDcat $ a \YD module,
and assume that the braiding of $V$ is of diagonal type.
Let $\cB (V)$ denote the corresponding Nichols algebra.
Our aim is to give an answer to the following question of
Andruskiewitsch stated in his survey \cite{inp-Andr02}.

\addvspace{\baselineskip}

\textsc{Question} 5.40. Given a braided vector space $V$ of
diagonal type and dimension 2, decide when $\cB (V)$ is finite
dimensional. If so, compute $\dim \cB (V)$, and give a
``nice'' presentation by generators and relations.

\addvspace{\baselineskip}

In the previous part we proved
finite dimensionality of several rank 2 Nichols algebras
with help of Gra\~nas differential operators
\cite{a-Grana00}, Kharchenkos theory \cite{a-Khar99} and
full binary trees. Here it is shown that any finite dimensional
rank 2 Nichols algebra of diagonal type is isomorphic to
one of the examples in \cite{p-Heck04a}.

The starting point of our classification is the following
observation. Under certain assumptions one can find elements of
$\cB (V)$ such that the subalgebra generated by them maps
surjectively onto another Nichols algebra. The latter is
called a subquotient of $\cB (V)$. The existence of
subquotients can be used in a very simple way to prove
infinite dimensionality of Nichols algebras. For example,
if there is a proper subquotient of $\cB (V)$ which is
isomorphic to $\cB (V)$ or if the subquotient is infinite
dimensional then $\cB (V)$ is infinite dimensional.

The method of subquotients presented here is applicable also for
Nichols algebras of higher rank. However our classification
method is effective only together with the constructive part
described in \cite{p-Heck04a}. The latter uses full binary trees
which have to be generalized or replaced if the rank of the 
Nichols algebra is greater than 2. It remains a challenging
problem to find an appropriate structure which carries
sufficiently many information needed to prove finiteness
results and to read off the algebra structure of $\cB (V)$.

The structure of the paper is as follows. First we formulate
general conditions for the existence of subquotients. In order
to allow an easier application some Corollaries are formulated.
Then in Theorem \ref{t-class} the list of all finite
dimensional rank 2 Nichols algebras is given, this time sorted
systematically by the entries of the braiding.
The main part of this paper is devoted to the description
of sufficiently many kinds of subquotients of $\cB (V)$.
In the last section the classification is splitted into
six special cases. In each
of them additional careful choices of subquotients are needed
in order to obtain valuable criterions for the entries of the
braiding.

We use the notation and conventions in \cite{p-Heck04a},
Section 2, which follow mainly \cite{a-AndrSchn98}.

\section{Nichols algebras and subquotients}
\label{s-Nicholssub}

Suppose that $k$ is a field of characteristic zero,
$G$ an abelian group, and $V\in {}^{kG}_{kG}\YDcat $
a finite dimensional \YD module with completely reducible
$kG$-action. Set $d:=\dim _kV$. Let $\lcoa :V\to kG\ot V$,
$\lact :kG \ot V\to V$, and $\sigma \in \End _k(V\ot V)$ denote the left
coaction of $kG$ on $V$, the left action of $kG$ on $V$, and the braiding
of $V$, respectively. Let $\cB (V)$ be the Nichols algebra generated by $V$.
The aim of this paper is to determine all $V\in {}^{kG}_{kG}\YDcat $
such that $\dim _kV=2$ and $\dim _k \cB (V)<\infty $.

Let $V^*$ denote the \YD module (left) dual to $V$.
Then the vector space $\cB (V^*)\ot kG$ together with the product
\begin{align*}
(f'\ot g')(f''\ot g''):=f'(g'\lact f'')\ot g'g'',\qquad
f',f''\in \cB (V^*),\,g',g''\in G,
\end{align*}
and the coproduct
\begin{align*}
\copr (f)&:=f\ot 1+\lcoa (f),& \copr (g):=g\ot g,& &f\in V^*,g\in G,
\end{align*}
is a Hopf algebra and will be denoted as usual by $\cB (V^*)\# kG$.
Recall the following lemma from \cite{p-Heck04a}.

\begin{lemma}\label{l-nichpaarung}
There exists a unique action $\paar{\cdot }{\cdot }$ of $\cB (V^*)\# kG$
on $\cB (V)$ satisfying
\begin{align*}
\paar{f}{v}=f(v), \quad \paar{g}{\rho }=g\lact \rho &
& &\text{for }f\in V^*,v\in V,g\in G,\rho \in \cB (V),\\
\paar{f_1f_2}{\rho }=\paar{f_1}{\paar{f_2}{\rho }}&
& &\text{for }f_1,f_2\in \cB (V^*)\# kG,\rho \in \cB (V),\\
\paar{f}{\rho \rho '}=\paar{f_{(1)}}{\rho }\paar{f_{(2)}}{\rho '}&
& &\text{for }f\in \cB (V^*)\# kG,\rho ,\rho '\in \cB (V).
\end{align*}
\end{lemma}

Suppose that $d\in \ndN $, $g_i\in G$ for $1\le i\le d$,
$q_{ij}\in k\setminus \{0\}$
for $i,j\in \{1,2,\ldots ,d\}$, and
$\{x_i\,|\,1\le i\le d\}$ is a basis of $V$ such that
\begin{align*}
\lcoa (x_i)&=g_i\ot x_i,&  g_i\lact x_j&=q_{ij}x_j.
\end{align*}
Such a basis always exists and is called \textit{a canonical basis} of $V$.
The algebra $\cB (V)$ admits an $\ndN _0^d$-grading such that $\deg x_i=
\Ndbasis _i$ where $\{\Ndbasis _i\,|\,1\le i\le d\}$ is a basis of the
$\ndN _0$-module $\ndN _0^d$. The corresponding $\ndZ $-grading, the so
called total grading, is defined by $\mathrm{totdeg}\,\rho :=\sum _{i=1}^d
n_i$ whenever $\rho \in \cB (V)$, $\deg \rho =\sum _{i=1}^dn_i\Ndbasis _i$.
The $\ndN _0^d$-grading and the left $kG$-coaction on $V$ induce
a group homomorphism $g:\ndZ ^d\to G$ and a bicharacter $\chi :\ndZ ^d\times
\ndZ ^d\to k$ given by
\begin{align*}
g(\Ndbasis _i)&:=g_i,& \chi (\Ndbasis _i,\Ndbasis _j):=q_{ij}.
\end{align*}
For notational convenience we will also write $g(x)$ and $\chi (x',x'')$
instead of $g(\deg x)$ and $\chi (\deg x',\deg x'')$ for homogeneous
elements $x,x',x''\in \cB (V)$.

Let $\{y_i\,|\,1\le i\le d\}$ denote the dual basis of $V^*$.
Then one gets
\begin{align}\label{eq-dualYD}
\lcoa (y_i)&=g_i^{-1}\ot y_i,& g_i\lact y_j&=q_{ij}^{-1}y_j,&
\sigma (y_i\ot y_j)&=q_{ij}y_j\ot y_i.
\end{align}
Thus for diagonal braidings the linear map $\iota :V\to V^*$,
$\iota (x_i):=y_i$ for $1\le i\le d$, extends to an algebra isomorphism
$\iota :\cB (V)\to \cB (V^*)$. 

The corollaries of the following lemma are our main tools in the
classification of Nichols algebras.

\begin{lemma}\label{l-surj}
Let $V$ and $W$ be finite dimensional \YD modules of diagonal type
over an abelian group $G$ and $H$, respectively. Set $d_V:=\dim _kV$ and
$d_W:=\dim _kW$.
Fix a canonical basis $\{w_i\,|\,1\le i\le d_W\}$ of $W$ and elements
$h_i\in H$, $1\le i\le d_W$, such that
$\lcoa (w_i)=h_i\ot w_i$. Further, let
$V_0\subset \cB (V)$ and $W_0\subset \cB (W)$ be \YD submodules of dimension
$n<\infty $.
Let $A$ and $B$ denote the subalgebras of $\cB (V)$ and $\cB (W)$ generated by
$V_0$ and $W_0$, respectively.
Choose canonical bases $\{v^0_i\,|\,1\le i\le n\}$ and
$\{w^0_i\,|\,1\le i\le n\}$ of $V_0$ and $W_0$, respectively.
Suppose that
\begin{itemize}
\item
$\paar{\iota (w_i)}{B}\subset B$ whenever $1\le i\le d_W$,
\item
there exists a $\ndZ $-grading $\deg _0$ of $\cB (V)$ such that
the elements $v^0_i$ are homogeneous with respect to $\deg _0$ and
$\deg _0(v^0_i)>0$ for $1\le i\le n$,
\item
there exist algebra automorphisms $\alpha _i\in \mathrm{Aut}_k(A)$,
and $k$-linear maps
$Y_i:A\to A$, $1\le i\le d_W$, and $\varphi _0: V_1\to B$ where
$V_1:=Y(V_0)$ and $Y$ denotes the unital $k$-subalgebra of
$\mathrm{End}_k(A)$ generated by $\{Y_i\,|\,1\le i\le d_W\}$, satisfying the
following properties.
\begin{itemize}
\item
For $1\le i\le d_W$ there exist $n_i\in \ndN $ such that
the maps $\alpha _i$ and $Y_i$ are homogeneous of degree
0 and $-n_i$ with respect to $\deg _0$, respectively,
\item
$Y_i(ab)=Y_i(a)b+\alpha _i(a)Y_i(b)$ for all $a,b\in A$
and $1\le i\le d_W$,
\item
$\varphi _0(v^0_j)=w^0_j$,
$\varphi _0(\alpha _i(v^0_j))=
h_i^{-1}\lact w^0_j$ for $1\le i\le d_W$ and $1\le j\le n$,
\item $\varphi _0(Y_i(a))=\paar{\iota (w_i)}{\varphi _0(a)}$
for all $a\in V_1$ and $1\le i\le d_W$.
\end{itemize}
\end{itemize}
Then there exists a unique surjective $k$-algebra homomorphism
$\varphi :A\to B$ such that $\varphi \restr{V_1}=\varphi _0$.
\end{lemma}

\begin{bew}
Since $\{v^0_i\,|\,1\le i\le n\}$ and $\{w^0_i\,|\,1\le i\le n\}$ generate $A$
and $B$, respectively, and $\varphi (v^0_i)=w^0_i$, surjectivity and
uniqueness of $\varphi $ is clear. To prove existence one has to show that
\begin{align*}
\tag{$*$}
a:=\sum _{m\in \ndN _0,1\le i_1,i_2,\ldots ,i_m\le n}&
\lambda _{m,i_1i_2\ldots i_m}v^0_{i_1}v^0_{i_2}\cdots v^0_{i_m}=0,\quad
\lambda _{m,i_1i_2\ldots i_m}\in k\\
&\Rightarrow
\sum _{m\in \ndN _0,1\le i_1,i_2,\ldots ,i_m\le n}
\lambda _{m,i_1i_2\ldots i_m}w^0_{i_1}w^0_{i_2}\cdots w^0_{i_m}=0.
\end{align*}
This can be done by induction over $\deg _0(a)$. Let $\cG _m$ denote the
set of homogeneous elements of $A$ of degree $m$ with respect to $\deg _0$.
Since $\deg _0(v^0_i)>0$ one has $\cG _0=k1$ and $\cG _m=\{0\}$ for $m<0$.
Thus the assertion ($*$) holds
for $a\in \cG _m$, $m\le 0$. Suppose now that $j\in \ndN _0$ and
($*$) holds for all $a\in \cG _m$, $m\le j$. Then for arbitrary
$a\in \cG _{j+1}$ we have to show that
$\paar{\iota (w_i)}{\tilde{a}}=0$ for $1\le i\le d_W$ where
$\tilde{a}:=\sum _{m\in \ndN _0,1\le i_1,i_2,\ldots ,i_m\le n}
\lambda _{m,i_1i_2\ldots i_m}w^0_{i_1}w^0_{i_2}\cdots w^0_{i_m}$.
One computes
\begin{align*}
\paar{\iota (w_i)}{&w^0_{i_1}w^0_{i_2}\cdots w^0_{i_m}}=
\sum _{l=0}^{m-1}(h^{-1}_i\lact w^0_{i_1}\ldots w^0_{i_l})
\paar{\iota (w_i)}{w^0_{i_{l+1}}}w^0_{i_{l+2}}\cdots w^0_{i_m}\\
=&\sum _{l=0}^{m-1}\varphi (\alpha _i(v^0_{i_1}\ldots v^0_{i_l}))
\varphi _0(Y_i(v^0_{i_{l+1}}))\varphi (v^0_{i_{l+2}}\ldots v^0_{i_m})
=\varphi (Y_i(v^0_{i_1}v^0_{i_2}\cdots v^0_{i_m})).
\end{align*}
Here the last two equations are valid because of the induction hypothesis
and the assumption on the degrees of $Y_i$ and $\alpha _i$. Now $a=0$ implies
$Y_i(a)=0$ for $1\le i\le d_W$ and hence $\paar{\iota (w_i)}{\tilde{a}}=0$.
\end{bew}

\begin{folg}\label{f-infinite}
In the setting of Lemma \ref{l-surj} suppose additionally that $V=W$
and $A$ is a proper subspace of $B$. Then $B$ and hence $\cB (V)$ are
infinite dimensional vector spaces.
\end{folg}

Let $\cB (V)^+$ denote the unique maximal ideal of $\cB (V)$.

\begin{folg}\label{f-subYD}
Let $V$ be a \YD module of diagonal type and $W\subset \cB (V)^+$
an $n$-dimensional \YD submodule where $n\in \ndN _0$.
Choose a canonical basis $\{w_i\,|\,1\le i\le n\}$ of $W$ and
let $A$ denote the subalgebra of $\cB (V)$ generated by $W$.
If there exist $\lambda _i\in k\setminus \{0\}$, $1\le i\le n$, such that
the restrictions of $\paar{\iota(w _i)}{\cdot }$ onto $A$ are
skew-primitive and $\paar{\iota (w _i)}{w_j}=\delta _{ij}\lambda _i$
then $\cB (W)$ is a quotient algebra of $A$. In particular,
$\dim _k \cB (W)=\infty $ implies that $\dim _k\cB (V)=\infty $.
\end{folg}

\begin{defin}
In the situation of Corollary \ref{f-subYD} we say that $\cB (W)$
\textit{is a subquotient of} $\cB (V)$ and write $\cB (W)\subq \cB (V)$.
If further $W\not= \emptyset $ and $W\not=V$ then we say that $\cB (W)$
\textit{is a proper subquotient of} $\cB (V)$ and write
$\cB (W)\psubq \cB (V)$.
\end{defin}

\begin{bew}[ of Corollary \ref{f-subYD}]
In Lemma \ref{l-surj} set $V_0=W_0:=W$ and $v^0_i:=w^0_i:=w_i$ for
$1\le i\le n$. For the grading $\deg _0$ take the grading
$\mathrm{totdeg}$ on $\cB (V)$. Further, set
$\alpha _i:=h_i^{-1}\lact (\cdot )$ and $Y_i:=\lambda _i^{-1}
\paar{\iota (w_i)}{\cdot }$. Then Lemma \ref{l-surj} gives a surjective
algebra homomorphism $\varphi :A\to \cB (W)$ extending the identity
$\varphi _0=\id :W\to W$.
\end{bew}

\section{The main result}
\label{s-main}

{}From now on consider only the case where $\dim _kV=2$. For $n\in \ndN $,
$n\ge 2$ let $R_n$ denote the subset of $k$ consisting of the primitive
$n^\mathrm{th}$ roots of unity.
The following theorem is the main result of this paper.

\begin{thm}\label{t-class}
Let $k$ be a field of characteristic zero and $G$ an abelian group.
Let $V\in {}_{kG}^{kG}\YDcat $ be a \YD module with $\dim _kV=2$ and completely
reducible $kG$-action. Assume that the Nichols algebra $\cB (V)$ is finite
dimensional. Then there exists a canonical basis $B$ of $V$ such that the
entries of the matrix $(q_{ij})_{i,j=1,2}$ of the braiding of $V$
with respect to $B$ satisfy
\begin{gather*}
q_{12}q_{21}=1 \text{ or }q_{12}q_{21}q_{22}=1 \text{ or }
q_{22}=-1\text{ or }q_{22}\in R_3.
\end{gather*}
More precisely, $B$ can be chosen such that one of the following holds.
\begin{enumerate}
\item
$q_{12}q_{21}=1$ and $q_{11},q_{22}\in \cup _{n=2}^\infty R_n$. 
\item
If $q_{12}q_{21}\not=1$, $q_{12}q_{21}q_{22}=1$:
\begin{itemize}
\item
$q_{11}q_{12}q_{21}=1$, $q_{12}q_{21}\in \cup _{n=2}^\infty R_n$.
\item
$q_{11}=-1$, $q_{12}q_{21}\in \cup _{n=3}^\infty R_n$.
\item
$q_{11}\in R_3$,
$q_{12}q_{21}\in \cup _{n=2}^\infty R_n$,
$q_{11}q_{12}q_{21}\not=1$.
\item
$q_{11}\in \cup _{n=4}^\infty R_n$,
$q_{12}q_{21}\in \{q_{11}^{-2},q_{11}^{-3}\}$.
\item
$q_{12}q_{21}\in R_8$, $q_{11}=(q_{12}q_{21})^2$.
\item
$q_{12}q_{21}\in R_{24}$, $q_{11}=(q_{12}q_{21})^6$.
\item
$q_{12}q_{21}\in R_{30}$, $q_{11}=(q_{12}q_{21})^{12}$.
\end{itemize}
\item
If $q_{12}q_{21}\not=1$, $q_{11}q_{12}q_{21}\not=1$,
$q_{12}q_{21}q_{22}\not=1$, $q_{22}=-1$, $q_{11}\in R_2\cup R_3$:
\begin{itemize}
\item
$q_{11}=-1$, $q_{12}q_{21}\in \cup _{n=3}^\infty R_n$.
\item
$q_{11}\in R_3$, $q_{12}q_{21}\in \{q_{11},-q_{11}\}$.
\item
$q_0:=q_{11}q_{12}q_{21}\in R_{12}$, $q_{11}=q_0^4$.
\item
$q_{12}q_{21}\in R_{12}$, $q_{11}=-(q_{12}q_{21})^2$.
\item
$q_{12}q_{21}\in R_9$, $q_{11}=(q_{12}q_{21})^{-3}$.
\item
$q_{12}q_{21}\in R_{24}$, $q_{11}=-(q_{12}q_{21})^4$.
\item
$q_{12}q_{21}\in R_{30}$, $q_{11}=-(q_{12}q_{21})^5$.
\end{itemize}
\item
If $q_{12}q_{21}\not=1$, $q_{11}q_{12}q_{21}\not=1$,
$q_{12}q_{21}q_{22}\not=1$, $q_{22}=-1$, $q_{11}\notin R_2\cup R_3$:
\begin{itemize}
\item
$q_{11}\in \cup _{n=5}^\infty R_n$, $q_{12}q_{21}=q_{11}^{-2}$.
\item
$q_{11}\in R_5\cup R_8\cup R_{12}\cup R_{14}\cup R_{20}$,
$q_{12}q_{21}=q_{11}^{-3}$.
\item
$q_{11}\in R_{10}\cup R_{18}$, $q_{12}q_{21}=q_{11}^{-4}$.
\item
$q_{11}\in R_{14}\cup R_{24}$, $q_{12}q_{21}=q_{11}^{-5}$.
\item
$q_{12}q_{21}\in R_8$, $q_{11}=(q_{12}q_{21})^{-2}$.
\item
$q_{12}q_{21}\in R_{12}$, $q_{11}=(q_{12}q_{21})^{-3}$.
\item
$q_{12}q_{21}\in R_{20}$, $q_{11}=(q_{12}q_{21})^{-4}$.
\item
$q_{12}q_{21}\in R_{30}$, $q_{11}=(q_{12}q_{21})^{-6}$.
\end{itemize}
\item
If $q_{12}q_{21}\not=1$, $q_{11}q_{12}q_{21}\not=1$,
$q_{12}q_{21}q_{22}\not=1$, $q_{11}\not=-1$, $q_{22}\in R_3$:
\begin{itemize}
\item
$q_0:=q_{11}q_{12}q_{21}\in R_{12}$, $q_{11}=q_0^4$, $q_{22}=-q_0^2$.
\item
$q_{12}q_{21}\in R_{12}$, $q_{11}=q_{22}=-(q_{12}q_{21})^2$.
\item
$q_{12}q_{21}\in R_{24}$, $q_{11}=(q_{12}q_{21})^{-6}$,
$q_{22}=(q_{12}q_{21})^{-8}$.
\item
$q_{11}\in R_{18}$, $q_{12}q_{21}=q_{11}^{-2}$, $q_{22}=-q_{11}^3$.
\item
$q_{11}\in R_{30}$, $q_{12}q_{21}=q_{11}^{-3}$, $q_{22}=-q_{11}^5$.
\end{itemize}
\end{enumerate}
\end{thm}

\begin{bem}
The list of Nichols algebras given in Theorem \ref{t-class} coincides
with the one of \cite[Theorem 4]{p-Heck04a}.
Here the examples are sorted systematically by their structure
constants whereas in \cite[Theorem 4]{p-Heck04a} they are listed by their
tree types. Therefore the Nichols algebras given in Theorem \ref{t-class}
are precisely all finite dimensional rank two Nichols algebras of diagonal
type over the field $k$.
\end{bem}

\section{Construction of subquotients}
\label{s-theory}

As in Section \ref{s-Nicholssub} let $\{x_1,x_2\}$ denote a canonical
basis of $V$. In this paper we want to classify
all Nichols algebras $\cB (V)$ satisfying
\begin{align*}
\tag{A0}
\dim _k\cB (V)<\infty .
\end{align*}
To do so we determine the necessary additional conditions for the structure
constants $q_{ij}$, $i,j\in \{1,2\}$.

\mysubs
The results in this subsection
essentially coincide with the assertion of Lemma 3.7 in \cite{inp-AndrSchn02}.

Note that $\paar{y_i}{x_i^m}=\sum _{j=0}^{m-1}q_{ii}^{-j}x_i^{m-1}$ and the
characteristic of $k$ is zero. Thus assumption (A0) implies that
\begin{align*}
\tag{A1}
q_{11},q_{22}\in \bigcup _{n=2}^\infty R_n.
\end{align*}
Set $x_{21}:=x_1x_2-q_{12}x_2x_1$ and $y_{21}:=\iota (x_{21})$. Then by
Equation (\ref{eq-dualYD}) and Lemma \ref{l-nichpaarung} one gets
\begin{gather*}
\paar{y_1}{x_{21}}=0, \quad \paar{y_2}{x_{21}}=(q_{21}^{-1}-q_{12})x_1,
\quad \paar{y_{21}}{x_{21}}=q_{21}^{-1}-q_{12},\\
\copr (y_{21})=y_{21}\ot 1+g_1^{-1}g_2^{-1}\ot y_{21}
+(1-q_{12}q_{21})y_1g_2^{-1}\ot y_2.
\end{gather*}
Now if (A0) holds then $x_{21}^m=0$ for some $m\in \ndN $. On the other hand
\begin{align*}
\paar{y_{21}}{x_{21}^m}&=\paar{y_{21}}{x_{21}}
\sum _{j=0}^{m-1}(q_{11}q_{12}q_{21}q_{22})^{-j}x_{21}^{m-1}.
\end{align*}
Thus (A0) implies that
\begin{align*}
\tag{A2}
q_{12}q_{21}=1\text{ or }
q_{11}q_{12}q_{21}q_{22}\in \bigcup _{n=2}^\infty R_n.
\end{align*}

Let $z_0:=x_2$, $z_{i+1}:=x_1z_i-q_{11}^iq_{12}z_ix_1$, and
$\hatz _i:=\iota (z_i)$ for $i\in \ndN _0$. Set $b_0:=1$,
$b_i:=\prod _{j=0}^{i-1}(1-q_{11}^jq_{12}q_{21})$ for $i\in \ndN $,
$\qnum{0}{p}:=0$, $\qnum{i}{p}:=\sum _{j=0}^{i-1}p^j$,
$\qfact{0}{p}:=1$, $\qfact{i}{p}:=\qnum{1}{p}\qnum{2}{p}\cdots \qnum{i}{p}$
for $i\in \ndN $ and $p\in k$, and
$\binom{i}{j}_p:=\qfact{i}{p}/(\qfact{j}{p}\qfact{i-j}{p})$ for
$i,j\in \ndN _0$, $j\le i$. Then one can prove by induction over $i$ that
\begin{gather}
\notag
\paar{y_1}{z_i}=0,\quad
\paar{y_2}{z_i}=q_{21}^{-i}b_ix_1^i,\quad
\paar{\hatz _j}{z_i}=\paar{y_1^jy_2}{z_i}=
\frac{q_{21}^{-i}b_i\qfact{i}{q_{11}^{-1}}}{\qfact{i-j}{q_{11}^{-1}}}
x_1^{i-j},\\
\label{eq-copr-hz}
\copr (\hatz _i)=\hatz _i\ot 1+\sum _{m=0}^i
\binom{i}{m}_{q_{11}}\frac{b_i}{b_m}y_1^{i-m}g_1^{-m}g_2^{-1}\ot \hatz _m
\end{gather}
for all $i,j\in \ndN _0$ with $j\le i$. In particular,
\begin{align}\label{eq-z=0}
z_i=0 \Leftrightarrow
b_i\qfact{i}{q_{11}}=0 \Leftrightarrow
\paar{\hatz _i}{z_i}=0.
\end{align}
Note that for all $i\in \ndN $ one has
\begin{align}\label{eq-di+10-di0}
\frac{\paar{\hatz _{i+1}}{z_{i+1}}}{\paar{\hatz _i}{z_i}}-
\frac{\paar{\hatz _i}{z_i}}{\paar{\hatz _{i-1}}{z_{i-1}}}=&
\chi (z_i,x_1)^{-1}-\chi (x_1,z_i).
\end{align}

Similarly let $u_0:=x_1$, $u_{i+1}:=u_ix_2-q_{12}q_{22}^ix_2u_i$, and
$\hatu _i:=\iota (u_i)$ for $i\in \ndN _0$. Set $c_0:=1$ and
$c_i:=\prod _{j=0}^{i-1}(1-q_{12}q_{21}q_{22}^j)$ for $i\in \ndN $.
Then one proves by induction over $i$ that
\begin{align*}
\paar{y_1}{u_i}=&\delta _{i0},&
\paar{y_2}{u_i}=&q_{21}^{-1}(1-q_{12}q_{21}q_{22}^{i-1})
\qnum{i}{q_{22}^{-1}}u_{i-1}
\end{align*}
for all $i\in \ndN _0$. In particular, $u_i=0$ if and only if
$c_i\qfact{i}{q_{22}}=0$.
Because of the symmetry of the conditions
for $u_i=0$ and $z_i=0$ one can choose the order of the basis vectors
$x_1$ and $x_2$ in such a manner that
\begin{align*}
\tag{A3}
\min\{i\in \ndN \,|\,u_i=0\}\le \min\{i\in \ndN \,|\,z_i=0\}
\end{align*}
whenever (A0) holds.
Thus assumptions (A0) and (A3) imply that either $z_1=0$ or $z_2=u_2=0$
or $z_2\not=0$, i.\,e.\
\begin{align*}
\tag{A4}
\text{$q_{12}q_{21}=1$ or $(1{-}q_{11}q_{12}q_{21})(1{+}q_{11})
=(1{-}q_{12}q_{21}q_{22})(1{+}q_{22})=0$ or $z_2\not=0$.}
\end{align*}

\mysubs
Assume that $z_i\not=0$ for some $i\in \ndN $. To shorten notation
we define $p_i:=\chi (z_i,z_i)^{-1}$.
By Lemma \ref{l-nichpaarung} and (\ref{eq-copr-hz}) one gets
$\paar{\hatz _i}{z_i^m}=\qnum{m}{p_i}
\paar{\hatz _i}{z_i}z_i^{m-1}$ for all $m>0$. Since $z_i^m=0$ for some
$m\in \ndN $ and $\paar{\hatz _i}{z_i}\not=0$ assumption (A0) implies
$p_i\not=1$.

Suppose that $z_{i+1}\not=0$ for some $i\in \ndN $.
If $p_i=-1$ then
\begin{align*}
\paar{\hatz _{i-1}}{z_i^2}&=q_{21}^{-i}b_i\qfact{i}{q_{11}^{-1}}(x_1z_i
+\chi (z_{i-1},z_i)^{-1}z_ix_1)=
q_{21}^{-i}b_i\qfact{i}{q_{11}^{-1}}z_{i+1},\\
\paar{\hatz _{i+1}\hatz _{i-1}}{z_i^{2m}}&=
\paar{\hatz _i}{z_i}\sum _{j=0}^{m-1}
\paar{\hatz _{i+1}}{\chi (z_{i-1},z_i)^{-2j}
z_i^{2j}z_{i+1}z_i^{2(m-1-j)}}\\
&=\paar{\hatz _i}{z_i}\paar{\hatz _{i+1}}{z_{i+1}}
\sum _{j=0}^{m-1}\chi (z_i,z_i)^{-4j}z_i^{2(m-1)}
\end{align*}
since $\paar{y_1}{z_i}=\paar{y_1}{z_{i+1}}=\paar{\hatz _{i+1}}{z_i}=0$.
As the characteristic of $k$ is zero this yields that $z_i^{2m}\not=0$
for all $m>0$. Therefore (A0) implies that
\begin{align*}
\tag{A5}
\begin{aligned}
i\in \ndN ,\quad p_i^{-1}(=\chi (z_i,z_i))=1 \quad &\Rightarrow \quad z_i=0,\\
i\in \ndN ,\quad p_i^{-1}=-1 \quad &\Rightarrow \quad z_{i+1}=0.
\end{aligned}
\end{align*}

\mysubs
\label{ss-zi2}
For $i\in \ndN _0$ set $z_{i,1}:=z_{i+1}z_i-\chi (z_{i+1},z_i)z_iz_{i+1}$,
$\hatz _{i,1}:=\iota (z_{i,1})$, and
\begin{align*}
d_{i,0}:=&q_{21}^{-1}(1-q_{11}^iq_{12}q_{21})\qnum{i+1}{q_{11}^{-1}}+
\chi (z_i,z_{i+1})^{-1}-\chi (z_{i+1},z_i).
\end{align*}
Then one obtains
\begin{align*}
\paar{\hatz _i}{z_{i,1}}=&\paar{\hatz _i}{z_{i+1}}z_i-
\chi (x_1,z_i)z_i\paar{\hatz _i}{z_{i+1}}\\
&+(\chi (z_i,z_{i+1})^{-1}-\chi (z_{i+1},z_i))
\paar{\hatz _i}{z_i}z_{i+1}=d_{i,0}\paar{\hatz _i}{z_i}z_{i+1}
\end{align*}
and hence $z_{i,1}=0$ implies $d_{i,0}z_{i+1}=0$.

Suppose now that $d_{i,0}z_{i+1}\not=0$ for some $i\in \ndN $ and
$p_i^2+p_i+1=0$. Then
\begin{align*}
\paar{\hatz _{i-1}}{z_i^3}=&\paar{\hatz _i}{z_i}(x_1z_i^2+
\chi (z_{i-1},z_i)^{-1}z_ix_1z_i+\chi (z_{i-1},z_i)^{-2}z_i^2x_1)\\
=&\paar{\hatz _i}{z_i}(z_{i+1}z_i-\chi (x_1,z_i)\chi (z_i,z_i)z_ix_1z_i
+\chi (z_{i-1},z_i)^{-2}z_i^2x_1)\\
=&\paar{\hatz _i}{z_i}(z_{i+1}z_i-\chi (z_{i+1},z_i)z_iz_{i+1})
=\paar{\hatz _i}{z_i}z_{i,1},
\end{align*}
\begin{align*}
\paar{\hatz _{i+1}\hatz _i\hatz _{i-1}}{z_i^{3m}}=&
\paar{\hatz _i}{z_i}\sum _{j=0}^{m-1}\chi (z_{i-1},z_i)^{-3j}
\paar{\hatz _{i+1}\hatz _i}{z_i^{3j}z_{i,1}z_i^{3(m-1-j)}}\\
=&d_{i,0}\paar{\hatz _i}{z_i}^2\sum _{j=0}^{m-1}\chi (z_iz_{i-1},z_i)^{-3j}
\paar{\hatz _{i+1}}{z_i^{3j}z_{i+1}z_i^{3(m-1-j)}}\\
=&md_{i,0}\paar{\hatz _i}{z_i}^2\paar{\hatz _{i+1}}{z_{i+1}}z_i^{3(m-1)}
\end{align*}
since $\paar{y_1}{z_i}=\paar{y_1}{z_{i,1}}=\paar{y_1}{z_{i+1}}
=\paar{\hatz _i}{z_i^3}=\paar{\hatz _{i+1}}{z_i}=0$
and $p_i^3=1$.
As the characteristic of $k$ is zero this yields that $z_i^{3m}\not=0$
for all $m>0$. Therefore using (A5) assumption (A0) implies that
\begin{align*}
\tag{A6}
i\in \ndN ,\quad p_i^3=1 \quad \Rightarrow
\quad d_{i,0}z_{i+1}=0.
\end{align*}

\mysubs
In this subsection assume only that (A5) holds.
Let $i\in \ndN $. If $z_{i+1}=0$ then set $w_i:=0$.
Otherwise by Equation (\ref{eq-z=0}) and (A5) one can define
\begin{align*}
w_i:=z_{i+1}z_{i-1}-\chi (z_{i+1},z_{i-1})z_{i-1}z_{i+1}
-\frac{\paar{\hatz _{i+1}}{z_{i+1}}}{\qnum{2}{p_i}\paar{\hatz _i}{z_i}}z_i^2.
\end{align*}
Set $\hatw _i:=\iota (w_i)$.
Our aim in this subsection is to prove the following lemma.

\begin{lemma}\label{l-w=0orinfty}
Assume that (A5) holds and let $i\in \ndN $.
If $\paar{\hatz _{m+1}\hatz _{m-1}}{w_m}=0$ for $1\le m\le i$ then
$w_m=0$ for $1\le m\le i$. Otherwise $\dim _k\cB (V)=\infty $.
\end{lemma}

Therefore (A0) implies the condition
\begin{align*}\tag{A7}
w_i=0 \quad \text{for all $i\in \ndN $.}
\end{align*}

In order to prove Lemma \ref{l-w=0orinfty} we can assume that
$z_{i+1}\not=0$ and hence
$z_2\not=0$. One obtains $\paar{y_1}{w_i}=0$ and
\begin{align*}
\paar{\hatz _{i-1}&}{w_i}=\frac{\paar{\hatz _{i+1}}{z_{i+1}}}{%
\qnum{2}{q_{11}^{-1}}}(x_1^2z_{i-1}-\chi (x_1^2,z_{i-1})z_{i-1}x_1^2)\\
&+\paar{\hatz _{i-1}}{z_{i-1}}(\chi (z_{i-1},z_{i+1})^{-1}-
\chi (z_{i+1},z_{i-1}))z_{i+1}\\
&-\frac{\paar{\hatz _{i+1}}{z_{i+1}}}{\qnum{2}{p_i}}
(x_1z_i+\chi (z_{i-1},z_i)^{-1}z_ix_1)\\
=&\bigg(\frac{\paar{\hatz _{i+1}}{z_{i+1}}}{\qnum{2}{q_{11}^{-1}}}
+\paar{\hatz _{i-1}}{z_{i-1}}\chi (z_{i-1},z_{i+1})^{-1}\\
&-\paar{\hatz _{i-1}}{z_{i-1}}\chi (z_{i+1},z_{i-1})
-\frac{\paar{\hatz _{i+1}}{z_{i+1}}}{\qnum{2}{p_i}}\bigg)z_{i+1}\\
=&(q_{11}p_i-1)\bigg(
\frac{q_{11}^{-1}\paar{\hatz _{i+1}}{z_{i+1}}}{%
\qnum{2}{q_{11}^{-1}}\qnum{2}{p_i}}
+q_{12}q_{21}^{-1}\paar{\hatz _{i-1}}{z_{i-1}}
(1+q_{11}^{-1}p_i^{-1})\bigg)z_{i+1}.
\end{align*}
In particular, one gets
\begin{gather*}
\paar{y_2}{w_1}=\frac{q_{21}^{-2}(1-q_{12}q_{21}q_{22})(1+q_{22}^{-1})
(1+q_{11}q_{12}^2q_{21}^2q_{22})}{1+q_{11}q_{12}q_{21}q_{22}}z_2,\\
\paar{\hatz _i}{w_i}=0,\qquad
\paar{\hatw _i}{w_i}=\paar{\hatz _{i+1}\hatz _{i-1}}{w_i}.
\end{gather*}
Further, let $A$ denote the
subalgebra of $\cB (V)$ generated by $w_i$ and $z_i$. In what follows
for $f\in \cB (V^*)\# kG$ let $f\restr{A}$ denote the map
$\paar{f}{\cdot }:A\to \cB (V)$. Then $\paar{y_1}{A}=0$
and hence $\hatz _j\restr{A}$ is skew-primitive for all $j\in \ndN _0$.
Since $\paar{\hatz _{i+1}}{A}=0$ one obtains
\begin{align*}
\copr &(\hatw _i)\restr{A\ot A}=\copr (\hatz _{i+1})(\hatz _{i-1}\ot 1
+g_1^{1-i}g_2^{-1}\ot \hatz _{i-1})\restr{A\ot A}\\
&-\frac{q_{21}^{-1}b_{i+1}\qnum{i+1}{q_{11}^{-1}}}{\qnum{2}{p_i}b_i}
\copr (\hatz _i)(\hatz _i\ot 1+g_1^{-i}g_2^{-1}\ot \hatz _i)\restr{A\ot A}.
\end{align*}
Using Equation (\ref{eq-copr-hz}), the definition of $\hatz _i$, and the fact
that $y_1\restr{A}=\hatz _m\restr{A}=0$ for $m>i$ one obtains
\begin{align*}
\copr (\hatz _{i+1})(\hatz _{i-1}\ot 1)&\restr{A\ot A}=\hatz _{i+1}
\hatz _{i-1}\restr{A}\ot 1\\
&+\qnum{i+1}{q_{11}}\frac{b_{i+1}}{b_i}\chi (z_i,z_{i-1})(\hatz _ig_1^{-i}
g_2^{-1}\ot \hatz _i)\restr{A\ot A},\\
\copr (\hatz _{i+1})(g_1^{1-i}g_2^{-1}\ot \hatz _{i-1})&\restr{A\ot A}=
(g_1^{-2i}g_2^{-2}\ot \hatz _{i+1}\hatz _{i-1})\restr{A\ot A},\\
\copr (\hatz _i)(\hatz _i\ot 1)&\restr{A\ot A}=(\hatz _i^2\ot 1
+\chi (z_i,z_i)\hatz _ig_1^{-i}g_2^{-1}\ot \hatz _i)\restr{A\ot A},\\
\copr (\hatz _i)(g_1^{-i}g_2^{-1}\ot \hatz _i)&\restr{A\ot A}=
(\hatz _ig_1^{-i}g_2^{-1}\ot \hatz _i
+g_1^{-2i}g_2^{-2}\ot \hatz _i^2)\restr{A\ot A}.
\end{align*}
These equations give that
\begin{align*}
\copr (\hatw _i)\restr{A\ot A}=(\hatw _i\ot 1+g_1^{-2i}g_2^{-2}\ot \hatw _i)
\restr{A\ot A}.
\end{align*}
Thus $\hatw _i\restr{A}$ is a skew-primitive endomorphism
of $A$.

\begin{bsp}
Assume that $g_2\in G$, $q\in k\setminus \{0\}$,
$g_1=g_2^2$, and the braiding of $V$ is given by the matrix
$\begin{pmatrix} q^4 & q^2\\ q^2 & q\end{pmatrix}$. Further, suppose that (A0)
holds.

1. By (A1) $q$ has to be a root of unity and $q^4\not=1$. Thus $z_1\not=0$.

2. If $q^5=1$ and $q\not=1$ then $z_2\not=0$ and $\chi (z_2,z_2)=1$
which contradicts (A5). Thus $u_2\not=0$.

3. If $q^4=-1$ then $u_2\not=0$ and $\chi (u_2,u_2)=1$. 
One gets $\paar{y_1y_2^2}{u_2^m}=m\paar{y_1y_2^2}{u_2}u_2^{m-1}$
which is a contradiction to (A0) and $\paar{y_1y_2^2}{u_2}\not=0$.
Thus $z_2\not=0$.

4. By (A5) with $i=1$ one has $q^9\not=-1$.

5. If $q^{13}=-1$ and $q\not=-1$ then $z_i\not=0$ for all $i<13$.
Further, one gets $\chi (z_6,z_6)=-1$ which is a contradiction to (A5).

6. Otherwise $\paar{y_2}{w_1}$ and $x_{21}$ are both nonzero. Set $W:=kw_1
+kx_{21}$. By Corollary \ref{f-subYD} one obtains $\cB (W)\psubq \cB (V)$
and the braiding of $\cB (W)$ with respect to the basis $\{w_1,x_{21}\}$
of $W$ is given by the matrix
$\begin{pmatrix} q^{36} & q^{18}\\ q^{18} & q^9\end{pmatrix}$.

Thus there are two possibilities.
\begin{itemize}
\item
There exists an infinite chain $\cdots \psubq \cB (V_n)
\psubq \cB (V_{n-1})\psubq \cdots \cB (V_1)\psubq \cB (V)$
of Nichols algebras where $V_i$,
$i\in \ndN $, are two dimensional \YD modules of diagonal type.
\item
There exists a finite chain
$\cB (V_n)\psubq \cB (V_{n-1})\psubq \cdots \cB (V_1)\psubq \cB (V_0)=\cB (V)$
of Nichols algebras where $V_i$,
$0\le i\le n$, are two dimensional \YD modules of diagonal type,
and $\cB (V_n)$ is infinite dimensional.
\end{itemize}
Therefore $\cB (V)$ is infinite dimensional.
\end{bsp}

\begin{bsp}
Assume that $g_1\in G$, $q\in k\setminus \{0\}$,
$g_2=g_1^3$, and the braiding of $V$ is given by the matrix
$\begin{pmatrix} q & q^3\\ q^3 & q^9\end{pmatrix}$.
Further, suppose that (A0) holds and $q^3\not=-1$.

1. By (A1) $q$ has to be a root of unity and $q^9\not=1$. Since $q^3\not=-1$
this yields $z_1\not=0$.

2. If $q\in R_7$ then $u_2\not=0$ and $\chi (u_2,u_2)=1$.
One gets $\paar{y_1y_2^2}{u_2^m}=m\paar{y_1y_2^2}{u_2}u_2^{m-1}$
for all $m>0$ which contradicts (A0). Hence by (\ref{eq-z=0})
one obtains that $z_2\not=0$.

3. If $q\in R_{18}$ then $z_3\not=0$ and $\chi (z_3,z_3)=1$
which contradicts (A5).

4. If $q^{15}=1$ and $q^3\not=1$ then $q\in R_5\cup R_{15}$.
In the first case $z_2\not=0$ and
$\chi (z_2,z_2)=1$ which is a contradiction to (A5). In the second case
one has $z_i=0$ if and only if $i\ge 10$. Further, $\chi (z_2,z_2)^3=q^{75}=1$
and $d_{2,0}=q^{-3}(1-q^8)(1+q^{-1}+q^{-2})\not=0$. This is a contradiction
to (A6).

5. If $q^{22}=-1$ then $q^2=-1$ or $q\in R_{44}$.
The first case is a contradiction to (A2). In the second case one
has $z_i=0$ if and only if $i\ge 39$. Further, $\chi (z_{19},z_{19})=1$
which is a contradiction to (A5).

6. Otherwise $\paar{y_2}{w_1}$ and $x_{21}$ are both nonzero. Set
$W:=kw_1+kx_{21}$.
By Corollary \ref{f-subYD} one obtains that $\cB (W)\psubq \cB (V)$ and the
braiding of $\cB (W)$ is given by the matrix
$\begin{pmatrix} q^{64} & q^{32}\\ q^{32} & q^{16}\end{pmatrix}$.
By the first example $\dim _k\cB (W)=\infty $.

Therefore $\dim _k\cB (V)<\infty $ implies $q^3=-1$.
\end{bsp}

We continue with the proof of Lemma \ref{l-w=0orinfty}. Suppose that
$\paar{\hatz _{i-1}}{w_i}\not=0$ for $i\in \ndN $. Then one can apply
Corollary \ref{f-subYD} with $W=kw_i+kz_i$ which gives $\cB (W)\psubq 
\cB (V)$. By Example 1 one gets $\dim _k\cB (V)=\infty $. This proves the
second part of the Lemma.

Assume that $\paar{\hatz _{m-1}}{w_m}=0$ for $1\le m\le i$.
We show by induction over $m$ that $w_m=0$ for $1\le m\le i$.
If $m=1$ then $\hatz _0=y_2$ and hence $w_1=0$. Now turn to the induction step.

We prove by induction on $-n$ that $\paar{\hatz _n}{w_{m+1}}=0$ for
$0\le n\le m$. Then the case $n=0$ proves that $w_{m+1}=0$.
By assumption $\paar{\hatz _m}{w_{m+1}}=0$. If $\paar{\hatz _{n+1}}{w_{m+1}}=0$
then $\paar{y_1}{\paar{\hatz _n}{w_{m+1}}}=0$ and hence there exists
$\lambda \in k$ such that $\paar{\hatz _n}{w_{m+1}}=\lambda z _{2m+2-n}$.
Thus it suffices to show that $\paar{\hatz _{2m+2-n}\hatz _n}{w_{m+1}}
(=\paar{y_1^{2m+2-n}y_2\hatz _n}{w_{m+1}})=0$.
Since $n+2\le 2m+2-n$ it suffices to check that
$\paar{y_1^{n+2}y_2\hatz _n}{w_{m+1}}(=\paar{\hatz _{n+2}\hatz _n}{w_{m+1}})
=0$. Since $n<m$ one has
$\hatw _{n+1}=0$ and hence the latter equation follows from the hypothesis
$\paar{\hatz _{n+1}}{w_{m+1}}=0$.

\mysubs
\begin{lemma}\label{l-zz=0}
Assume that (A7) is satisfied. Let $i,j\in \ndN _0$ such that $z_i\not=0$ and
$j<i-1$. Then
$z_iz_j-\chi (z_i,z_j)z_jz_i\in \Lin _k\{z_mz_{i+j-m}\,|\,j<m<i\}$.
\end{lemma}

\begin{bew}
We use induction over $i-j$. If $i=j+2$ then the claim follows from
$w_{j+1}=0$. To prove the induction step note that
\begin{align*}
z_{i+1}z_j=&(x_1z_i-\chi (x_1,z_i)z_ix_1)z_j=x_1\chi (z_i,z_j)z_jz_i \\
&\text{ $+$ terms in
$x_1\Lin _k\{z_mz_{i+j-m}\,|\,j<m<i\}\oplus
kz_iz_{j+1}\oplus kz_iz_jx_1$}\\
=&\chi (z_{i+1},z_j)z_jz_{i+1}
\text{ $+$ terms in
$\Lin _k\{z_mz_{i+j+1-m}\,|\,j<m<i+1\}$}\\
& \text{ $+$ terms in
$\Lin _k\{z_mz_{i+j-m}x_1\,|\,j\le m\le i\}$}.
\end{align*}
Since $\paar{y_1}{z_{i+1}z_j}=0$ this gives the assertion.
\end{bew}

\begin{folg}\label{f-zi1=0}
If (A7) holds then  $z_{i,1}=0$ if and only if $d_{i,0}z_{i+1}=0$.
\end{folg}

\begin{bew}
Since $\paar{y_1}{z_{i,1}}=0$ one has $z_{i,1}=0$ if and only if
$\paar{\hatz _{j_1}\hatz _{j_2}}{z_{i,1}}=0$ whenever $j_1,j_2\in \ndN _0$
and $j_1+j_2=2i+1$. By Lemma \ref{l-zz=0} this is equivalent to
the fact that $\paar{\hatz _{j_1}\hatz _{j_2}}{z_{i,1}}=0$
whenever $j_1,j_2\in \ndN _0$, $j_1\le j_2+1$, and $j_1+j_2=2i+1$.
Since $\paar{\hatz _{i+1}}{z_{i,1}}=0$ this is the same as
$\paar{\hatz _{i+1}\hatz _i}{z_{i,1}}=0$.
\end{bew}

\begin{folg}\label{f-zisq=0}
If (A7) holds then  $z_i^2=0$ if and only if $\qnum{2}{p_i}z_i=0$.
\end{folg}

\begin{bew}
Analogous to the proof of Corollary \ref{f-zi1=0}.
\end{bew}

\mysubs
In this subsection assume that (A5)--(A7) hold.
Let $i\in \ndN $. If $z_{i,1}=0$ then set $s_i:=0$. Otherwise
$d_{i,0}z_{i+1}\not=0$ by Corollary \ref{f-zi1=0},
$p_i^3\not=1$ by (A6), and $p_i\not=-1$ by (A5) and since
$z_{i+1}\not=0$. Thus one can define
\begin{align*}
s_i:=&z_{i,1}z_{i-1}-\chi (z_{i,1},z_{i-1})z_{i-1}z_{i,1}
-\frac{d_{i,0}\paar{\hatz _{i+1}}{z_{i+1}}}{\qfact{3}{p_i}
\paar{\hatz _i}{z_i}}z_i^3.
\end{align*}
For $j\in \ndN $ let
$\cI (y_1,\hatz _j)$ denote the left ideal of $\cB (V^*)$ generated by
$y_1$ and $\hatz _j$. Then using (\ref{eq-copr-hz}) direct
computation shows that
\begin{gather}
\label{eq-copr-hz1}
\begin{aligned}
\copr (\hatz _{j,1})=&\hatz _{j,1}\ot 1 +g(z_{j,1})^{-1}\ot \hatz _{j,1}
+\chi (z_j,z_{j+1})d_{j,0}\hatz _{j+1}g(z_j)^{-1}\ot \hatz _j\\
&\text{$+$ terms in $(\cI (y_1,\hatz _{j+2})\# kG)\ot \cB (V^*)$}
\end{aligned}
\end{gather}
for all $j\in \ndN _0$.
In this subsection the following lemma will be proved.

\begin{lemma}\label{l-s=0orinfty}
Assume that (A5)--(A7) hold and let $i\in \ndN $.
If $\paar{\hatz _{i,1}\hatz _{i-1}}{s_i}=0$ then $s_i=0$.
Otherwise $\qnum{3}{-p_i}=0$ or $\dim _k\cB (V)=\infty $.
\end{lemma}

Let $i\in \ndN $. Without loss of generality assume that
$z_{i,1}\not=0$. By (A7) one has
$w_m=0$ for $m\in \ndN $. Further, $\paar{y_1}{s_i}=0$ and
\begin{align*}
\paar{\hatz _i}{s_i}=&\paar{\hatz _i}{z_{i,1}}z_{i-1}-
\chi (z_{i+1},z_{i-1})z_{i-1}\paar{\hatz _i}{z_{i,1}}
-d_{i,0}\paar{\hatz _{i+1}}{z_{i+1}}\qnum{2}{p_i}^{-1}z_i^2=0.
\end{align*}
By Lemma \ref{l-zz=0} one has $s_i=0$ if and only if $\paar{\hatz _{j_1}
\hatz _{j_2}\hatz _{j_3}}{s_i}=0$ whenever $j_1,j_2,j_3\in \ndN _0$,
$j_1+j_2+j_3=3i$, and $j_1\le j_2+1\le j_3+2$. Since $\paar{\hatz _i}{s_i}=0$
this yields that $s_i=0$ is equivalent to $\paar{\hatz _{i+1}\hatz _i
\hatz _{i-1}}{s_i}=0$.

Since $\paar{\hatz _i}{s_i}=0$ and $w_i=0$
one has $\paar{\hatz _{i,1}\hatz _{i-1}}{s_i}
=\paar{\hatz _{i+1}\hatz _i\hatz _{i-1}}{s_i}=\paar{\hatz _{i+1}
\hatz _{i-1,1}}{s_i}$. Moreover, Equation $\hatw _i=0$ implies that
\begin{align}
\notag
y_1\hatz _{i-1,1}=&(\hatz _{i+1}+\chi (x_1,z_i)\hatz _iy_1)\hatz _{i-1}
-\chi (z_i,z_{i-1})(\hatz _i+\chi (x_1,z_{i-1})\hatz _{i-1}y_1)\hatz _i\\
\notag
=&\left(\frac{\paar{\hatz _{i+1}}{z_{i+1}}}{\qnum{2}{p_i}\paar{\hatz _i}{z_i}}
+\chi (x_1,z_i)-\chi (z_i,z_{i-1})\right)\hatz _i^2
+\chi (x_1,z_{i-1,1})\hatz _{i-1,1}y_1\\
\label{eq-y1zi1}
=&d_{i-1,0}\qnum{2}{p_i}^{-1}\hatz _i^2
+\chi (x_1,z_{i-1,1})\hatz _{i-1,1}y_1
\end{align}
by (\ref{eq-di+10-di0}). Clearly one gets $\paar{\hatz _{i-1,1}}{z_{i,1}}
\in kx_1^2$. Further, $z_{i+1}\not=0$ implies that $z_2\not=0$ and
$\qnum{2}{q_{11}^{-1}}\not=0$, and hence Equations $\paar{y_1}{z_{i,1}}=0$
and (\ref{eq-y1zi1}) give
\begin{align*}
\paar{\hatz _{i-1,1}}{z_{i,1}}=&\qnum{2}{q_{11}^{-1}}^{-1}\paar{y_1
\hatz _{i-1,1}}{z_{i,1}}x_1
=\frac{d_{i-1,0}d_{i,0}\paar{\hatz _i}{z_i}\paar{\hatz _{i+1}}{z_{i+1}}}{
\qnum{2}{p_i}\qnum{2}{q_{11}^{-1}}}x_1^2.
\end{align*}
Using (\ref{eq-di+10-di0}) and Equations $\paar{\hatz _i}{z_{i-1,1}}=0$,
$\paar{\hatz _{i+1}}{z_jx_1^{i+1-j}}=0$ for $j\le i$, one gets
\begin{align*}
\paar{\hatz _{i,1}\hatz _{i-1}}{s_i}=&\paarb{\hatz _{i+1}\hatz _{i-1,1}}{
z_{i,1}z_{i-1}-\chi (z_{i,1},z_{i-1})z_{i-1}z_{i,1}-
\frac{d_{i,0}\paar{\hatz _{i+1}}{z_{i+1}}}{
\qfact{3}{p_i}\paar{\hatz _i}{z_i}}z_i^3}\\
=&\paarb{\hatz _{i+1}}{\paar{\hatz _{i-1,1}}{z_{i,1}}z_{i-1}+
d_{i-1,0}\chi (z_{i-1},z_{i+1})^{-1}\paar{\hatz _i}{z_{i,1}}
\paar{\hatz _{i-1}}{z_{i-1}}\\
&\quad -\frac{d_{i,0}\paar{\hatz _{i+1}}{z_{i+1}}}{\qfact{3}{p_i}
\paar{\hatz _i}{z_i}}
d_{i-1,0}\paar{\hatz _i}{z_i}\paar{\hatz _{i-1}}{z_i^2}}\\
=&d_{i-1,0}d_{i,0}\paar{\hatz _i}{z_i}\paar{\hatz _{i+1}}{z_{i+1}}^2\bigg(
\frac{1}{\qnum{2}{p_i}\qnum{2}{q_{11}^{-1}}}\\
&+\chi (z_{i-1},z_{i+1})^{-1}\frac{\paar{\hatz _{i-1}}{z_{i-1}}}{%
\paar{\hatz _{i+1}}{z_{i+1}}}-\frac{1}{\qfact{3}{p_i}}\bigg).
\end{align*}
To prove the second statement of the lemma it suffices to show that
$\iota (s_i)\restr{A}$ is skew-primitive where $A$ denotes the
subalgebra of $\cB (V)$ generated by $z_i$ and $s_i$. Indeed, in this case
Corollary \ref{f-subYD} with $W=kz_i+ks_i$ and Example 2 imply the assertion
since $p_i\not=-1$ and $\hatz _i\restr{A}$ is skew-primitive.

First note that
\begin{align*}
s_i=&z_{i+1}z_iz_{i-1}-\chi (z_{i+1},z_i)z_iz_{i+1}z_{i-1}
-\chi (z_{i,1},z_{i-1})z_{i-1}z_{i,1}-\frac{d_{i,0}
\paar{\hatz _{i+1}}{z_{i+1}}}{\qfact{3}{p_i}\paar{\hatz _i}{z_i}}z_i^3\\
=&z_{i+1}z_{i-1,1}+\chi (z_i,z_{i-1})z_{i+1}z_{i-1}z_i
-\chi (z_{i+1},z_{i-1,1})z_iz_{i-1}z_{i+1}\\
&-\chi (z_{i+1},z_i)
\frac{\paar{\hatz _{i+1}}{z_{i+1}}}{\qnum{2}{p_i}\paar{\hatz _i}{z_i}}z_i^3
-\chi (z_{i,1},z_{i-1})z_{i-1}z_{i,1}-\frac{d_{i,0}
\paar{\hatz _{i+1}}{z_{i+1}}}{\qfact{3}{p_i}\paar{\hatz _i}{z_i}}z_i^3\\
=&z_{i+1}z_{i-1,1}-\chi (z_{i+1},z_{i-1,1})z_{i-1,1}z_{i+1}\\
&+\left(\chi (z_i,z_{i-1})-\chi (z_{i+1},z_i)
-\frac{d_{i,0}}{\qnum{3}{p_i}}\right)\frac{\paar{\hatz _{i+1}}{z_{i+1}}}{%
\qnum{2}{p_i}\paar{\hatz _i}{z_i}}z_i^3\\
=&z_{i+1}z_{i-1,1}-\chi (z_{i+1},z_{i-1,1})z_{i-1,1}z_{i+1}
-\frac{d_{i-1,0}\paar{\hatz _{i+1}}{z_{i+1}}}{\qfact{3}{p_i}
\paar{\hatz _i}{z_i}}z_i^3.
\end{align*}
Now one has $\paar{\hatz _{i+1}}{A}=0$ and $\paar{\hatz _i}{A}\subset A$
and one computes
\begin{align*}
\copr (\hatz _{i+1}&\hatz _{i-1,1})\restr{A\ot A}=\copr (\hatz _{i+1})
(\hatz _{i-1,1}\ot 1 +g(z_{i-1,1})^{-1}\ot \hatz _{i-1,1}\\
&+\chi (z_{i-1},z_i)d_{i-1,0}\hatz _ig(z_{i-1})^{-1}\ot \hatz _{i-1})
\restr{A\ot A}\\
=&\bigg(\hatz _{i+1}\hatz _{i-1,1}\ot 1+\chi (z_i,z_{i-1,1})\qnum{i+1}{q_{11}}
\frac{b_{i+1}}{b_i}y_1\hatz _{i-1,1}g(z_i)^{-1}\ot \hatz _i\\
&+g(s_i)^{-1}\ot \hatz _{i+1}\hatz _{i-1,1}
+\chi (z_i,z_i)^2d_{i-1,0}\hatz _ig(z_i)^{-2}\ot \hatz _{i+1}\hatz _{i-1}\Big)
\restr{A\ot A}\\
\overseteqnum{eq-y1zi1}{=}&\Big(\hatz _{i+1}\hatz _{i-1,1}\ot 1
+g(s_i)^{-1}\ot \hatz _{i+1}\hatz _{i-1,1}\\
&+\chi (z_i,z_i)^2\frac{d_{i-1,0}\paar{\hatz _{i+1}}{z_{i+1}}}{\qnum{2}{p_i}
\paar{\hatz _i}{z_i}}(\hatz _i^2 g(z_i)^{-1}\ot \hatz _i
+\hatz _ig(z_i)^{-2}\ot \hatz _i^2)\bigg)\restr{A\ot A},\\
\copr (\hatz _i^3&)\restr{A\ot A}=\copr (\hatz _i)^3\restr{A\ot A}=
(\hatz _i^3\ot 1+\qnum{3}{p_i^{-1}}\hatz _i^2g(z_i)^{-1}\ot \hatz _i\\
&+\qnum{3}{p_i^{-1}}\hatz _ig(z_i)^{-2}\ot \hatz _i^2
+g(z_i)^{-3}\ot \hatz _i^3)\restr{A\ot A}.
\end{align*}
This proves that $\iota (s_i)\restr{A}$ is skew-primitive.

\mysubs
In this subsection assume that (A5) and (A7) hold.
Set
\begin{align*}
z_{i,2}:=z_{i+1}z_{i,1}-\chi (z_{i+1},z_{i,1})z_{i,1}z_{i+1},\quad
i\in \ndN _0.
\end{align*}
Using Equation (\ref{eq-copr-hz}) and $\hatw _{i+1}=0$
one obtains immediately the formulas
\begin{align}\label{eq-paarzi2}
\paar{\hatz _{i+1}}{z_{i,2}}&=0, &
\paar{\hatz _{i+2}\hatz _i}{z_{i,2}}&=0, &
\paar{\hatz _{i+1}\hatz _i}{z_{i,2}}&\in kz_{i+1}.
\end{align}
For $i\in \ndN _0$ set
\begin{align*}
d_{i,1}:=q_{21}^{-1}(1-q_{11}^iq_{12}q_{21})
\qnum{i+1}{q_{11}^{-1}}+\qnum{2}{p_{i+1}}
(\chi (z_i,z_{i+1})^{-1}-\chi (z_{i+1},z_iz_{i+1})).
\end{align*}
Then using the formulas for
$\copr (\hatz _i)$, $\copr (\hatz _{i,1})$ and Equation (\ref{eq-y1zi1})
one obtains that
\begin{gather}
\label{eq-copr-hz2}
\begin{aligned}
\copr (\hatz _{i,2})-\hatz _{i,2}\ot 1-g(z_{i,2})^{-1}\ot \hatz _{i,2}
-\chi (z_{i,1},z_{i+1})d_{i,1}\hatz _{i+1}g(z_{i,1})^{-1}\ot \hatz _{i,1}&\\
-\chi (z_i,z_{i+1})^2d_{i,0}d_{i,1}\qnum{2}{p_{i+1}}^{-1}
\hatz _{i+1}^2g(z_i)^{-1}\ot \hatz _i&
\end{aligned}
\end{gather}
is an element of $(\cI (y_1,\hatz _{i+2},\hatz _{i+2}\hatz _{i+1})\# kG)
\ot \cB (V^*)$.
Here $\cI (y_1,\hatz _{i+2},\hatz _{i+2}\hatz _{i+1})$ denotes the left ideal
of $\cB (V^*)$ generated by $y_1$, $\hatz _{i+2}$, and
$\hatz _{i+2}\hatz _{i+1}$.
The expression (\ref{eq-copr-hz2}) does not make sense if
$\qnum{2}{p_{i+1}}=0$. In this case $z_{i+2}=0$ by (A5) and $\hatz _{i+1}^2=0$
by Corollary \ref{f-zisq=0}. Thus $y_1\hatz _{i,1}=\chi (x_1,z_{i,1})
\hatz _{i,1}y_1$ and the formula for $\copr (\hatz _{i,2})$ takes the form
(\ref{eq-copr-hz2}) without the last summand.

Now we prove that for all $i\in \ndN _0$ one has
\begin{align}
\label{eq-normzi2}
\paar{\hatz _{i+1}\hatz _i}{z_{i,2}}=d_{i,0}d_{i,1}\paar{\hatz _i}{z_i}
\paar{\hatz _{i+1}}{z_{i+1}}z_{i+1}.
\end{align}
Recall from Subsection \ref{ss-zi2} that $\paar{\hatz _i}{z_{i,1}}=
d_{i,0}\paar{\hatz _i}{z_i}z_{i+1}$. Therefore one gets
\begin{align*}
\paar{\hatz _{i+1}\hatz _i}{z_{i,2}}
\overseteqnum{eq-paarzi2}{=}&
\paar{\hatz _{i,1}}{z_{i+1}z_{i,1}-\chi (z_{i+1},z_{i,1})z_{i,1}z_{i+1}}\\
\overseteqnum{eq-copr-hz1}{=}&\paar{\hatz _{i,1}}{z_{i,1}}(\chi(z _{i,1},
z _{i+1})^{-1}-\chi (z_{i+1},z_{i,1}))z_{i+1}\\
&+d_{i,0}\paar{\hatz _{i+1}}{z_{i+1}}\paar{\hatz _i}{z_{i,1}}\\
=&d_{i,0}\paar{\hatz _i}{z_i}\paar{\hatz _{i+1}}{z_{i+1}}
(\chi (z_{i,1},z_{i+1})^{-1}-\chi (z_{i+1},z_{i,1})+d_{i,0})z_{i+1}
\end{align*}
which implies (\ref{eq-normzi2}).

\begin{lemma}
Assume that (A5) and (A7) hold. Let $i\in \ndN _0$ such that $z_{i+1}\not=0$.
Then one has
\begin{align*}
\paar{\hatz _i}{z_{i,2}}=\begin{cases}
d_{i,0}d_{i,1}\qnum{2}{p_{i+1}}^{-1}\paar{\hatz _i}{z_i}z_{i+1}^2
&\text{if $p_{i+1}\not=-1$,}\\
0 &\text{if $p_{i+1}=-1$.}
\end{cases}
\end{align*}
\end{lemma}

\begin{bew}
By (\ref{eq-paarzi2}) one has $\paar{\hatz _{i+1}}{z_{i,2}}=0$.
Therefore Lemma \ref{l-zz=0} implies that
$\paar{\hatz _i}{z_{i,2}}\in
\Lin _k\{z_{i+1-n}z_{i+1+n}\,|\,0\le n\le i+1\}$.
Note that $z_{i+1+n}=0$ if and only if $\paar{\hatz _{i+1}}{z_{i+1+n}}=0$.
Further, $\paar{\hatz _{i+1}\hatz _i}{z_{i,2}}\in kz_{i+1}$ by
(\ref{eq-normzi2}). This together with Equation
\begin{align*}
\paar{\hatz _{i+1}}{z_{i+1-n}z_{i+1+n}}&=
\chi (z_{i+1},z_{i+1-n})^{-1}z_{i+1-n}\paar{\hatz _{i+1}}{z_{i+1+n}},\quad
n\in \ndN ,
\end{align*}
and the fact that the latter is an element of
$(k\setminus \{0\})z_{i+1}x_1^n$, implies that
$\paar{\hatz _i}{z_{i,2}}\in kz_{i+1}^2$. Now apply
$\paar{\hatz _{i+1}}{\cdot }$ and use (\ref{eq-normzi2}) and
Corollary \ref{f-zisq=0}.
\end{bew}

\begin{folg}\label{f-zi2=0}
If (A5) and (A7) hold then $z_{i,2}=0$ if and only if
$d_{i,0}d_{i,1}\qnum{2}{p_{i+1}}z_{i+1}=0$.
\end{folg}

\begin{bew}
Similar to the proof of Corollary \ref{f-zi1=0}.
\end{bew}

\mysubs
\label{ss-4.7}
Suppose that (A5) and (A7) hold, $\chi (z_{i,1},z_{i,1})=-1$
for some $i\in \ndN _0$,
and $z_{i,2}\not=0$. Then one has
$\paar{\hatz _{i,1}}{z_{i,2}}\in kz_{i+1}$ and
$\paar{\hatz _{i,1}}{z_{i,2}}\not=0$ by Corollary
\ref{f-zi2=0} and Equations (\ref{eq-paarzi2}) and (\ref{eq-normzi2}).
{}From this one gets $\paar{\hatz _{i,2}}{z_{i,2}}\not=0$. The latter
implies that $z_{i,1}^{2m}\not=0$ for all $m>0$. Indeed, one computes
$\paar{\hatz _i}{z_{i,1}^2}=d_{i,0}\paar{\hatz _i}{z_i}z_{i,2}$
and $\paar{\hatz _{i,2}\hatz _i}{z_{i,1}^{2m}}=m\paar{\hatz _{i,2}
\hatz _i}{z_{i,1}^2}z_{i,1}^{2m-2}$.

\mysubs
To obtain sufficiently many subquotients we will need the following
expressions. Assume again that (A5)--(A7) hold and
fix $i\in \ndN $ such that $\chi (z_{i,1},z_{i,1})\not=-1$
(see also \ref{ss-4.7}) and $z_{i,2}\not=0$. Set
\begin{align*}
t_i:=&z_{i,2}z_i-\chi (z_{i,2},z_i)z_iz_{i,2}
-\frac{d_{i,1}}{1+\chi (z_{i,1},z_{i,1})^{-1}}z_{i,1}^2
\end{align*}
and $\hatt _i:=\iota (t_i)$. Our aim is to apply Corollary \ref{f-subYD}
with $W=kt_i+kz_{i,1}$.

One has $\paar{\hatz _{i+1}}{t_i}=0$.
Let $A$ denote the subalgebra of $\cB (V)$ generated
by the elements $t _i$ and $z _{i,1}$.
Equations (\ref{eq-copr-hz1}) and (\ref{eq-normzi2}) yield
that $\hatz _{i,1}\restr{A}$ is skew-primitive and
$\paar{\hatz _{i,1}}{t_i}=0$. Therefore
$\paar{\hatz _{i,1}}{A}\subset A$.
One gets $\paar{\hatz _{i,2}}{t_i}=0$
as well and hence $\hatz _{i,2}\restr{A}=0$ by (\ref{eq-copr-hz2}).

Note that $\paar{\hatz _{i,1}}{A}\subset A$ implies
$\hatz _{i+m}\hatz _{i+1}\hatz _i\restr{A}=
\hatz _{i+m}\hatz _{i,1}\restr{A}=0$ for $m\ge 1$. Therefore
using the formulas for $\copr (\hatz _i)$, $\copr (\hatz _{i,1})$, and
$\copr (\hatz _{i,2})$ one gets
\begin{align*}
\copr (\hatt _i)&\restr{A\ot A}=\copr (\hatz _{i,2})
(\hatz _i\ot 1+g(z_i)^{-1}\ot \hatz _i)\restr{A\ot A}\\
&-\frac{d_{i,1}}{1+\chi (z_{i,1},z_{i,1})^{-1}}
\copr (\hatz_{i,1})(\hatz _{i,1}\ot 1
+g(z_{i,1})^{-1}\ot \hatz _{i,1})\restr{A\ot A}.
\end{align*}
For the summands of this expression one computes
\begin{gather*}
\copr (\hatz _{i,2})(\hatz _i\ot 1)\restr{A\ot A}=
(\hatz _{i,2}\hatz _i\ot 1+\chi (z_{i,1},z_{i,1})d_{i,1}
\hatz _{i,1}g(z_{i,1})^{-1}\ot \hatz _{i,1})\restr{A\ot A},\\
\copr (\hatz _{i,2})(g(z_i)^{-1}\ot \hatz _i)\restr{A\ot A}=
(g(t_i)^{-1}\ot \hatz _{i,2}\hatz _i)\restr{A\ot A},\\
\copr (\hatz _{i,1})(\hatz _{i,1}\ot 1)\restr{A\ot A}=
(\hatz _{i,1}^2\ot 1+\chi (z_{i,1},z_{i,1})\hatz _{i,1}g(z_{i,1})^{-1}
\ot \hatz _{i,1})\restr{A\ot A},\\
\copr (\hatz _{i,1})(g(z_{i,1})^{-1}\ot \hatz _{i,1})
\restr{A\ot A}=(g(z_{i,1})^{-2}\ot \hatz _{i,1}^2
+\hatz _{i,1}g(z_{i,1})^{-1}\ot \hatz _{i,1})\restr{A\ot A}.
\end{gather*}
Therefore
$\copr (\hatt _i)\restr{A\ot A}=(\hatt _i\ot 1
+g(t_i)^{-1}\ot \hatt _i)\restr{A\ot A}$.

To apply Corollary \ref{f-subYD} with $W=kt_i+kz_{i,1}$
one has to check the relation
$\paar{\hatt _i}{t _i}\not=0$. Since
$\paar{\hatz _{i,1}}{t_i}=0$ one gets
$\paar{\hatt _i}{t_i}=\paar{\hatz _{i,2}\hatz _i}{t_i}$.
Recall that there exist $\mu _m\in k$ such that
$\paar{\hatz _i}{z_{i,m}}=\mu _mz_{i+1}^m$ for $m\in \{1,2\}$.
Therefore
\begin{align*}
\paar{\hatz _i}{t_i}\in kz_{i,2}+kz_{i,1}z_{i+1}+kz_iz_{i+1}^2.
\end{align*}
Since $\paar{\hatz _{i,1}}{t_i}=0$ one obtains
$\paar{\hatz _i}{t_i}\in kz_{i,2}$. Note that $z_{i,2}\not=0$
implies that $\qnum{2}{p_{i+1}}\not=0$
by Corollary \ref{f-zi2=0}.
Thus using the definition of $t_i$ one gets
\begin{align*}
\paar{\hatz _i}{t_i}=&\bigg(\mu _2+(\chi (z_i,z_{i,2})^{-1}
-\chi (z_{i,2},z_i))\paar{\hatz _i}{z_i}
-\frac{d_{i,1}\mu _1}{1+\chi (z_{i,1},z_{i,1})^{-1}}\bigg)z_{i,2}\\
=&
\paar{\hatz _i}{z_i}
\bigg(\frac{d_{i,0}d_{i,1}}{\qnum{2}{p_{i+1}}}
-\frac{d_{i,0}d_{i,1}}{1+\chi (z_{i,1},z_{i,1})^{-1}}
+\chi (z_i,z_{i,2})^{-1}-\chi (z_{i,2},z_i)\bigg)z_{i,2}.
\end{align*}

\begin{lemma}\label{l-t=0orinfty}
Suppose that (A5)--(A7) hold. Let $i\in \ndN $ such that $z_{i,2}\not=0$
and either $\chi (z_{i,1},z_{i,1})=-1$
or $\chi (z_{i,1},z_{i,1})\not=-1$, $\paar{\hatz _i}{t_i}\not=0$. Then
$\dim _k\cB (V)=\infty $.
\end{lemma}

\begin{bew}
If $\chi (z_{i,1},z_{i,1})=-1$ then $\dim _k\cB (V)=\infty $ by
Subsection \ref{ss-4.7}. Otherwise $\paar{\hatz _i}{t_i}\not=0$
implies that $\paar{\hatt _i}{t_i}\not=0$ and one can apply
Corollary \ref{f-subYD} with $W=kt_i+kz_{i,1}$.
Then by Example 1 one gets $\dim _k\cB (V)=\infty $.
\end{bew}

\section{The classification}

The structures worked out in the previous sections are sufficient
to perform the proposed classification of finite dimensional Nichols
algebras. Thus in this section we generally assume that $G$ is an abelian
group, $V\in {}^{kG}_{kG}\YDcat $ is a two-dimensional
\YD module of diagonal type
and $\cB (V)$ is the corresponding Nichols algebra such that (A0) and (A3)
hold. Suppose first that in (A4) we have $z_2=0$. Then by (A1) and (A2)
$\cB (V)$ appears in Theorem \ref{t-class}(1, 2.1, 2.2, 3.1). More precisely,
if $q_{11}q_{12}q_{21}=1$ and $q_{22}=-1$ then first one has to exchange
the variables $x_1$ and $x_2$. Thus to prove Theorem \ref{t-class}
it remains to consider the case when
$z_2\not=0$. Set $a:=\min \{i\in \ndN \,|\,z_{i+2}=0\}$.
By (A7) $w_a$ has to be zero.
The relations $z_{a+2}=0$, $z_{a+1}\not=0$, and Equation
(\ref{eq-z=0}) imply that
$(1-q_{11}^{a+1}q_{12}q_{21})(1-q_{11}^{a+2})=0$. If
$q_{11}^{a+1}q_{12}q_{21}=1$ then
\begin{align*}
\paar{\hatz _{a-1}}{w_a}=&\frac{(q_{11}^{a+1}q_{22}^{-1}-1)
\paar{\hatz _{a-1}}{z_{a-1}}}{1+q_{11}^aq_{22}^{-1}}
q_{11}^{-1}q_{21}^{-2}(1+q_{22}^{-1})(1+q_{11}^{-2a-1}q_{22})z_{a+1}
\end{align*}
and if $q_{11}^{a+2}=1$ then
\begin{align*}
\paar{\hatz _{a-1}}{w_a}=&(q_{11}^{-3}(q_{12}q_{21})^{-a}q_{22}^{-1}-1)
\frac{q_{11}^2q_{21}^{-2}\paar{\hatz _{a-1}}{z_{a-1}}}{1+
q_{11}^{-4}(q_{12}q_{21})^{-a}q_{22}^{-1}}\times \\
&\qquad (q_{11}(q_{12}q_{21})^{a+1}q_{22}+1)
(1+q_{11}^{-6}(q_{12}q_{21})^{1-a}q_{22}^{-1})z_{a+1}.
\end{align*}
Therefore $w_a=0$ allows only the following six possibilities.
\begin{itemize}
\item
$q_{22}=q_{11}^{a+1}$, $q_{12}q_{21}=q_{11}^{-a-1}$.
\item
$q_{22}=-1$, $q_{12}q_{21}=q_{11}^{-a-1}$.
\item
$q_{22}=-q_{11}^{2a+1}$, $q_{12}q_{21}=q_{11}^{-a-1}$.
\item
$q_{11}^{a+2}=1$, $q_{22}=q_{11}^{-3}(q_{12}q_{21})^{-a}$.
\item
$q_{11}^{a+2}=1$, $q_{22}=-q_{11}^{-6}(q_{12}q_{21})^{-a+1}$.
\item
$q_{11}^{a+2}=1$, $q_{22}=-q_{11}^{-1}(q_{12}q_{21})^{-a-1}$.
\end{itemize}
Note that by the definition of $a$ we may also assume that
\begin{align*}
\tag{A8}
\text{$q_{11}^m\not=1$ for all $m\in \ndN $, $m<a+2$.}
\end{align*}
In the following subsections we will analyze each case separately.

\mysubs
\label{ss-Serie1}
Consider the case
$q_{22}=q_{11}^{a+1}$, $q_{12}q_{21}=q_{11}^{-a-1}$, $a\ge 1$.
Then the braiding is of Cartan type (cf.~\cite[Sect.~4]{inp-AndrSchn02})
and hence it is known for many
values of $q_{11}$ that $\dim _k\cB (V)<\infty $ if and only if $a\le 2$
(see e.\,g.~Theorem~4.6 in \cite{inp-AndrSchn02}).
Note that by (A1) and (A8) one has $q_{11}\in \cup _{n=a+2}^\infty R_n$.
The settings for $a=1$ and $a=2$ appear in Theorem \ref{t-class}(2.4)
up to the case $a=1$, $q_{11}\in R_3$, which is a special case of
Theorem \ref{t-class}(2.3).
We prove $\dim _k\cB (V)=\infty $ for $a\ge 3$.

Since $\chi (z_2,z_2)=q_{11}^{4-2(a+1)+a+1}=q_{11}^{3-a}$ one gets a
contradiction to (A5) if $a=3$.

Let $a\ge 4$. Set $n:=2$, $v^0_1:=z_{a-1}$, $v^0_2:=z_a$, $w^0_1:=z_2$,
$w^0_2:=x_1$, $\deg _0(x_1):=-2$, $\deg _0(x_2):=2a+1$, and
\begin{gather*}
\alpha _1:=g_1^{-a}g_2^{-1}\lact (\cdot ),\qquad
\alpha _2:=g_1^{a+1}g_2\lact (\cdot ),\\
Y_1:=\frac{\paar{\hatz _a}{\cdot }}{\paar{\hatz _a}{z_a}},\quad
Y_2(\rho ):=\frac{\paar{\hatz _2}{z_2}\paar{\hatz _a}{z_a}}{%
\paar{\hatz _{a+1}}{z_{a+1}}}(z_{a+1}\rho -(g_1^{a+1}g_2\lact \rho )
z_{a+1})
\end{gather*}
for $\rho \in \cB (V)$. One gets the formulas
\begin{align*}
\chi (z_a,v^0_1)^{-1}=&q_{11}^{-2}q_{12}^{-1}=\chi (x_1,w^0_1)^{-1},&
\chi (z_a,v^0_2)^{-1}=&q_{11}^{-1}=\chi (x_1,w^0_2)^{-1},\\
\chi (z_{a+1},v^0_1)=&q_{11}^{-a-1}q_{21}^{-2}=\chi (x_2,w^0_1)^{-1},&
\chi (z_{a+1},v^0_2)=&q_{21}^{-1}=\chi (x_2,w^0_2)^{-1}.
\end{align*}
Further, $d_{a,0}=0$ and hence $z_{a,1}=0$ by Corollary \ref{f-zi1=0}.
Since $w_a=0$ as well the smallest subspace of $\cB (V)$ containing
$kv^0_1$ and $kv^0_2$ and stable under the action of $Y_1$ and $Y_2$
is $V_1:=kz_{a-1}+kz_a+kz_a^2$. Now one can check step by step that all
assumptions of Corollary \ref{f-infinite} are fulfilled
where one has to set $n_1:=1$, $n_2:=1$, and
$\varphi _0(z_{a-1}):=z_2$, $\varphi _0(z_a):=x_1$,
$\varphi _0(z_a^2):=x_1^2$.
This gives a contradiction to (A0).

\mysubs
\label{ss-Serie2}
Assume now that $q_{12}q_{21}=q_{11}^{-a-1}$ and $q_{22}=-1$.
Then $u_2=0$, $p_i=-q_{11}^{i(a+1-i)}$ for all $i$, and
\begin{align*}
d_{1,0}&=q_{21}^{-1}(1-q_{11}^{-a}+q_{11}^{-2a})(1-q_{11}^{a-1}),\\
d_{2,0}&=q_{21}^{-1}(1-q_{11}^{a-2})(q_{11}^{a-2}+1-q_{11}^{1-a}-q_{11}^{-a}
+q_{11}^{1-2a}+q_{11}^{3-3a}).
\end{align*}
If $a=1$ and $q_{11}\in R_4$ then $\cB (V)$ appears in Theorem
\ref{t-class}(2.4). Otherwise $a=1$ and (A8) imply that $q_{11}^2\not=1$
and $\cB (V)$ appears in Theorem \ref{t-class}(3.2, 4.1).

Further, if $a\ge 2$ then $\paar{y_2}{w_1}=0$ and
\begin{align*}
\paar{\hatz _1}{w_2}&=q_{21}^{-3}b_1\frac{q_{11}^{-a}(1+q_{11}^{2a-1})}{
1+q_{11}^{a-1}}(1-q_{11}^{-a}+q_{11}^{-2a})(1-q_{11}^{a-2})z_3.
\end{align*}
If $a=2$ then one gets
$d_{1,1}=q_{21}^{-1}(1-q_{11}^{-2})(1+q_{11}^3)(1+q_{11}^{-4})$
and $\qnum{2}{p_2}=1-q_{11}^2$.
Thus by (A8) and Corollary \ref{f-zi2=0} one has
$z_{1,2}=0$ if and only if $q_{11}\in R_6\cup R_8\cup R_{12}$.
The case $q_{11}\in R_6$ was already considered in Subsection
\ref{ss-Serie1}.
If $q_{11}\in R_8\cup R_{12}$ then $\cB (V)$ appears in Theorem
\ref{t-class}(4.2).
Otherwise by Lemma \ref{l-t=0orinfty}
$\paar{\hatz _1}{t_1}$ has to be zero. One computes
\begin{align*}
\frac{\paar{\hatz _1}{t_1}}{\paar{\hatz _1}{z_1}}=&
\frac{q_{11}^{-10}q_{21}^{-2}(1-q_{11}^5)(1+q_{11}^7)\qnum{5}{-q_{11}^2}}{%
1-q_{11}^3+q_{11}^6}z_{1,2}
\end{align*}
which yields $q_{11}\in R_5\cup R_{14}\cup R_{20}$. In these cases
$\cB (V)$ appears in Theorem \ref{t-class}(4.2).

Let $a\ge 3$. Then by (A8) and Equation $\paar{\hatz _1}{w_2}=0$ one obtains
that $(q_{11}^{2a-1}+1)(q_{11}^{2a}-q_{11}^a+1)=0$.
This means that if $a=3$ then $q_{11}\in R_{10}\cup R_{18}$ and
if $a=4$ then $q_{11}\in R_{14}\cup R_{24}$. In
these cases $\cB (V)$ appears in Theorem \ref{t-class}(4.3) and
Theorem \ref{t-class}(4.4), respectively.

Suppose now that $a\ge 5$. Then $w_3=0$ implies that
\begin{align*}
0=&\frac{\paar{\hatz _2}{w_3}}{\paar{\hatz _2}{z_2}}=
\frac{-q_{21}^{-2}(q_{11}^{3a-5}+1)(q_{11}^{a-3}-1)}{\qnum{3}{q_{11}^{a-2}}}
\times \\
&(q_{11}^{-2}+q_{11}^{1-a}+q_{11}^{-a}-q_{11}^{2-2a}-q_{11}^{1-2a}
-q_{11}^{-2a}+q_{11}^{2-3a}+q_{11}^{1-3a}+q_{11}^{4-4a})z_4.
\end{align*}
For $q_{11}^{2a-1}=-1$ and $\qnum{3}{-q_{11}^a}=0$ one obtains
\begin{align*}
\frac{\paar{\hatz _2}{w_3}}{\paar{\hatz _2}{z_2}}=&
\frac{-q_{21}^{-2}(1-q_{11}^{a-4})(q_{11}^{a-3}-1)}{\qnum{3}{q_{11}^{a-2}}}
q_{11}^{-2}\qnum{5}{q_{11}}z_4,\\
\frac{\paar{\hatz _2}{w_3}}{\paar{\hatz _2}{z_2}}=&
\frac{-q_{21}^{-2}(1-q_{11}^{-5})(q_{11}^{a-3}-1)}{\qnum{3}{q_{11}^{a-2}}}
(1+q_{11}^{-2})(1-q_{11}^{4-a})z_4,
\end{align*}
respectively. In both cases (A8) implies that $\paar{\hatz _2}{w_3}\not=0$
which is a contradiction to (A7).

\mysubs
\label{ss-Serie3}
In this subsection we consider the case $q_{12}q_{21}=q_{11}^{-a-1}$,
$q_{22}=-q_{11}^{2a+1}$. Note that by (A3) one has $u_i=0$
for some $i\in \{2,3,\ldots ,a+2\}$. If $u_2=0$ then $q_{22}=-1$
or $q_{12}q_{21}q_{22}=1$ and hence this case is already covered by the
previous subsections.

If $a=1$ then $u_3$ has to be zero as well. This gives $q_{11}^2=-1$ or
$q_{11}^6-q_{11}^3+1=0$. However $q_{11}^2=-1$ yields
$\chi (z_1,z_1)=1$ which is a contradiction to (A5).
If $q_{11}\in R_{18}$ then $\cB (V)$ appears in Theorem
\ref{t-class}(5.4).

For $a=2$ one obtains
\begin{align*}
d_{1,0}&=q_{11}q_{21}^{-1}\qnum{3}{q_{11}^{-1}}^2(1-q_{11}^{-1})\not= 0
\qquad \text{by (A8)},\\
d_{1,1}&=q_{11}^{-7}q_{21}^{-1}(1-q_{11})(1+q_{11}^2)(
q_{11}^8+q_{11}^7-q_{11}^5-q_{11}^4-q_{11}^3+q_{11}+1).
\end{align*}
Further, $u_2\not=0$ implies that $1\not=q_{12}q_{21}q_{22}=-q_{11}^2$
and (A5) and $z_3\not=0$ yield $p_2\not=-1$.
Therefore from Corollary \ref{f-zi2=0} we obtain that $z_{1,2}=0$
if and only if $q_{11}\in R_{30}$.
In this case $\cB (V)$ appears in Theorem \ref{t-class}(5.5).
Otherwise Lemma \ref{l-t=0orinfty} gives that
$\paar{\hatz _1}{t_1}$ has to be zero. Therefore
\begin{align*}
&\frac{d_{1,0}d_{1,1}}{1-q_{11}^{-3}}
-\frac{d_{1,0}d_{1,1}}{1+q_{11}^{-11}}
-q_{21}^{-2}q_{11}^{-11}+q_{21}^{-2}q_{11}^5=q_{11}^{-6}q_{21}^{-2}
\frac{(q_{11}^2+1)(1-q_{11}^{-5})}{q_{11}^{11}+1}\times \\
&\quad (q_{11}^8+1)
(q_{11}^{12}-q_{11}^{10}-q_{11}^9+q_{11}^7+q_{11}^6+q_{11}^5-q_{11}^3-q_{11}^2
+1)=0.
\end{align*}
If $q_{11}^5=1$ then $q_{22}=-1$. This case was considered in Subsection
\ref{ss-Serie2}. Further, as mentioned above $u_2\not=0$ implies
$q_{11}^2\not=-1$.

Recall that $u_2\not=0$, $u_3=0$ yields
$(q_{11}^7-1)\qnum{3}{-q_{11}^5}=0$ and $u_3\not=0$, $u_4=0$
implies that $(q_{11}^{12}+1)(1+q_{11}^{10})=0$.
Further, we assumed that $q_{11}\notin R_{30}$. Thus
Equations $\paar{\hatz _1}{t_1}=0$ and
$u_4=0$ lead to a contradiction.

Let $a\ge 3$. Then
\begin{align*}
\paar{\hatz _1}{w_2}=&\frac{-q_{11}^{-1}q_{21}^{-3}(1+q_{11}^{-2})b_1}{%
1-q_{11}^{-1}}(1-q_{11}^{2-a})(1-q_{11}^{-a-1})z_3.
\end{align*}
Thus Equation $w_2=0$ gives a contradiction to (A8).

\mysubs
\label{ss-Serie4}
In the fourth case we suppose that $q_{11}^{a+2}=1$ and
$q_{22}=q_{11}^{-3}(q_{12}q_{21})^{-a}$. By (A8) it can be
assumed that $q_{11}\in R_{a+2}$.
Further, the case $(q_{12}q_{21})^{a+2}=1$ can be excluded as well. Indeed,
if $(q_{12}q_{21})^{a+2}=1$ then $q_{11}^iq_{12}q_{21}=1$
for some $i\in \ndN _0$, $i\le a+1$. But $z_{a+1}\not=0$
implies that $i<a+1$ is not possible and the case $i=a+1$
was already considered in Subsections \ref{ss-Serie1}--\ref{ss-Serie3}.

If $a=1$ then $q_{11}\in R_3$, $q_{12}q_{21}q_{22}=1$, and
$(q_{12}q_{21})^3\not=1$, and $\cB (V)$ appears in Theorem
\ref{t-class}(2.3).

For $a\ge 2$ one gets
\begin{align*}
\frac{\paar{\hatz _{a-2}}{w_{a-1}}}{\paar{\hatz _{a-2}}{z_{a-2}}}=&
\frac{q_{11}^2q_{21}^{-2}(q_{11}^{-5}q_{12}q_{21}-1)}{1+
q_{11}^{-6}q_{12}q_{21}}(q_{11}^2+1)(q_{11}+1)\qnum{3}{-q_{11}^{-4}q_{12}
q_{21}}z_a.
\end{align*}
Assume that $a=2$. Then $q_{11}^2=-1$,
$q_{22}=q_{11}(q_{12}q_{21})^{-2}$, and $w_1=0$. Further, by (A3)
$u_4$ has to be zero. If $u_2=0$ then
$(q_{22}+1)(q_{12}q_{21}q_{22}-1)=0$. If $q_{12}q_{21}q_{22}=1$ then
$q_{11}=q_{12}q_{21}$ which was excluded.
If $q_{22}=-1$ then $q_{12}q_{21}\in R_8$,
$q_{11}=(q_{12}q_{21})^{-2}$, and $\cB (V)$ appears in
Theorem \ref{t-class}(4.5).

If $a=2$, $u_2\not=0$, and $u_3=0$ then
$(q_{12}q_{21}q_{22}^2-1)\qnum{3}{q_{22}}=0$.
Since $(q_{12}q_{21})^4\not=1$ there are two possibilities:
\begin{itemize}
\item $q_{22}\in R_{12}$, $q_{11}=q_{22}^{-3}$,
$q_{12}q_{21}=q_{22}^{-2}$. Then $\chi (z_1,z_1)=q_{22}^{-4}$
and $d_{1,0}=q_{21}^{-1}(1-q_{22}^{-2})$ which is a
contradiction to (A6).
\item $q_{12}q_{21}\in R_{24}$, $q_{11}=(q_{12}q_{21})^{-6}$,
$q_{22}=(q_{12}q_{21})^{-8}$. Then $\cB (V)$ appears in
Theorem \ref{t-class}(5.3).
\end{itemize}
If $a=2$, $u_3\not=0$, and $u_4=0$ then
$(q_{22}^2+1)(q_{12}q_{21}q_{22}^3-1)=0$.
Again there are two possibilities:
\begin{itemize}
\item
$q_{11}^2=q_{22}^2=-1$. Then $q_{22}=q_{11}(q_{12}q_{21})^{-2}$
implies $(q_{12}q_{21})^4=1$. This is a contradiction to the
assumption at the beginning of this subsection.
\item
$q_{22}^{10}=-1$, $q_{12}q_{21}=q_{22}^{-3}$, $q_{11}=q_{22}^{-5}$.
The change of the role of $x_1$ and $x_2$ leads to an algebra
which was already considered in Subsection \ref{ss-Serie3}.
\end{itemize}

Suppose that $a\ge 3$. Then we must have
$\paar{\hatz _{a-2}}{w_{a-1}}=0$ by (A7).
If $q_{12}q_{21}=q_{11}^5$ then
$(q_{12}q_{21})^{a+2}=1$ which was excluded. Further, (A8) gives
$q_{11}^4\not=1$. Thus one obtains
$q_{11}^{-8}(q_{12}q_{21})^2-q_{11}^{-4}q_{12}q_{21}+1=0$
and hence $q_{11}^{-12}q_{12}^3q_{21}^3=-1$.
Then Equation
\begin{align*}
0=\frac{\paar{\hatz _{a-3}}{w_{a-2}}}{\paar{\hatz _{a-3}}{z_{a-3}}}=
\frac{q_{12}^2(q_{11}^{-12}q_{12}^2q_{21}^2-1)}{1+
q_{11}^{-13}q_{12}^2q_{21}^2}\qnum{3}{-q_{11}^{-1}}\qnum{5}{q_{11}^{-1}}
(q_{11}^{-7}q_{12}q_{21}-1)z_{a-1}
\end{align*}
together with $(q_{12}q_{21})^{a+2}\not=1$ and
$q_{11}^{-12}q_{12}^2q_{21}^2=-(q_{12}q_{21})^{-1}$ implies that
$\qnum{3}{-q_{11}}\qnum{5}{q_{11}}(q_{12}q_{21}+1)=0$.

If $\qnum{3}{-q_{11}}=0$ then $q_{11}^3=-1$. Since
$q_{11}\in R_{a+2}$ one gets $a=4$.
On the other hand, $(q_{11}^{-4}q_{12}q_{21})^3=-1$ implies
$(q_{12}q_{21})^3=-1$ which contradicts the assumption
$(q_{12}q_{21})^{a+2}\not=1$.

If $q_{12}q_{21}=-1$ then
$\qnum{3}{-q_{11}^{-4}q_{12}q_{21}}=0$ implies that
$q_{11}\in R_3\cup R_6\cup R_{12}$. Since
$a\ge 3$ this is again a contradiction to
$(q_{12}q_{21})^{a+2}\not=1$.

In the remaining case one has $a=3$, $\qnum{5}{q_{11}}=0$, and
$-q_{11}q_{12}q_{21}\in R_3$. This is
equivalent to $q_{12}q_{21}\in R_{30}$,
$q_{11}=(q_{12}q_{21})^{-6}$, which then yields
that $q_{22}=-1$. This example appears in
Theorem \ref{t-class}(4.8).

\mysubs
\label{ss-Serie5}
In this subsection assume that $q_{11}\in R_{a+2}$
and $q_{22}=-q_{11}^{-6}(q_{12}q_{21})^{1-a}$
where $a\in \ndN $. Further, as in Subsection
\ref{ss-Serie4} one can exclude the case $(q_{12}q_{21})^{a+2}=1$.
One has $p_{a+1}=-q_{11}^5(q_{12}q_{21})^{-2}$,
\begin{align*}
\chi (z_{a-1},z_{a-1})=&-q_{11}^3,&
\chi (z_{a},z_{a-1})=&-q_{12},&
\chi (z_{a},z_{a+1})=&-q_{11}^{-4}q_{12}q^2_{21},
\\
\chi (z_{a-1},z_{a})=&-q_{21},&
\chi (z_{a},z_{a})=&-q_{11}^{-2}q_{12}q_{21},&
\chi (z_{a+1},z_{a})=&-q_{11}^{-4}q_{12}^2q_{21}.
\end{align*}
\begin{align*}
d_{a-1,0}=&q_{21}^{-1}\qnum{3}{q_{11}}(q_{11}^{-2}q_{12}q_{21}-1),\\
d_{a,0}=&q_{11}q_{21}^{-1}(q_{11}^3(q_{12}q_{21})^{-1}+1)
(q_{11}^{-5}(q_{12}q_{21})^2-1),\\
d_{a,1}=&-q_{11}q_{21}^{-1}(q_{11}^{-10}(q_{12}q_{21})^4-q_{11}^{-5}
(q_{12}q_{21})^2+1)(1-q_{11}^8(q_{12}q_{21})^{-3}),\\
\frac{\paar{\hatz _a}{t_a}}{\paar{\hatz _a}{z_a}}=&
\frac{q_{11}^{-15}q_{12}^{-3}q_{21}^{-5}((q_{12}q_{21})^4+q_{11}^{10})}{%
1+q_{11}^{-5}(q_{12}q_{21})^2}
(q_{12}q_{21}-q_{11}^2)((q_{12}q_{21})^5+q_{11}^{13})z_{a,2}.
\end{align*}
Suppose that $a=1$. Then $q_{11}\in R_3$ and $q_{22}=-1$. Consider
the equation $d_{1,0}=0$. The case $q_{12}q_{21}=-1$ was part of
the previous subsection. On the other hand, if
$(q_{12}q_{21})^2=q_{11}^2$ then $\cB (V)$ appears in Theorem
\ref{t-class}(3.2).
Otherwise $d_{1,0}\not=0$ and $\qnum{2}{p_2}\not=0$.
Since $q_{11}^2(q_{12}q_{21})^4-q_{11}(q_{12}q_{21})^2
+1=((q_{12}q_{21})^2+1)(q_{11}^2(q_{12}q_{21})^2+1)$
there are three possibilities for $d_{1,1}=0$.
First if $q_{12}^2q_{21}^2=-1$ then set $q_0:=q_{11}q_{12}q_{21}$. One gets
$\qnum{3}{-q_0^2}=0$, $q_{11}=q_0^4$, $q_{12}q_{21}=-q_0^3$, $q_{22}=-1$.
This example appears in Theorem \ref{t-class}(3.3).
Next if $(q_{12}q_{21})^2=-q_{11}$ then $q_{12}q_{21}\in R_{12}$
and $q_{22}=-1$. In this case $\cB (V)$ appears in Theorem
\ref{t-class}(3.4).
Finally, if $q_{11}=(q_{12}q_{21})^{-3}$ then $q_{12}q_{21}\in R_9$,
$q_{22}=-1$, and $\cB (V)$ appears in Theorem \ref{t-class}(3.5).
If nothing of these equations is true then by Corollary \ref{f-zi2=0}
$z_{1,2}\not=0$ and hence
$\paar{\hatz _1}{t_1}$ has to vanish. Since $(q_{12}q_{21})^3\not=1$
one has $q_{12}q_{21}\not=q_{11}^2$. Since $\chi (z_1,z_1)\not=1$ by (A5) one
obtains $q_{11}q_{12}q_{21}\not=-1$. Therefore there are two remaining cases.
If $\qnum{3}{-(q_{12}q_{21})^4}=0$ and $q_{11}=-(q_{12}q_{21})^4$ then
$\cB (V)$ appears in Theorem \ref{t-class}(3.6).
On the other hand, if $q_{12}q_{21}\in R_{30}$ and
$q_{11}=-(q_{12}q_{21})^5$ then
$\cB (V)$ appears in Theorem \ref{t-class}(3.7).

If $a=2$ then $q_{11}^2=-1$ and $q_{22}=(q_{12}q_{21})^{-1}$.
Lemma \ref{l-s=0orinfty} implies that $\qnum{3}{-q_{12}q_{21}}=0$ or
$\paar{\hatz _{2,1}\hatz _1}{s_2}=0$. One computes
\begin{align*}
\paar{\hatz _{2,1}\hatz _1}{s_2}=&\frac{
d_{1,0}d_{2,0}\paar{\hatz _2}{z_2}\paar{\hatz _3}{z_3}^2
q_{12}q_{21}\qnum{3}{q_{11}^{-1}q_{12}^{-2}q_{21}^{-2}}}{(1+q_{11}^{-1})
(1+q_{12}^{-1}q_{21}^{-1})\qnum{3}{q_{12}^{-1}q_{21}^{-1}}
(1-q_{11}q_{12}q_{21})}.
\end{align*}
Since $(q_{12}q_{21})^4\not=1$ one has $d_{1,0}\not=0$. 
Relations $d_{2,0}=0$ and $(q_{12}q_{21})^4\not=1$ imply that
$(q_{12}q_{21})^4=-1$, $q_{11}=(q_{12}q_{21})^2$.
In this case $\cB (V)$ appears in Theorem \ref{t-class}(2.5).
If $\qnum{3}{q_{11}^{-1}q_{12}^{-2}q_{21}^{-2}}=0$ then
$q_{12}q_{21}\in R_{24}$, $q_{11}=(q_{12}q_{21})^6$,
and $\cB (V)$ appears in Theorem \ref{t-class}(2.6).
On the other hand, if
$\qnum{3}{-q_{12}q_{21}}=0$ then set $q_0:=q_{11}q_{12}q_{21}$. Now one has
$q_{11}=q_0^3$, $q_{12}q_{21}=-q_0^4$, $\qnum{3}{-q_0^2}=0$,
$d_{2,0}=-q_{21}^{-1}q_0(1+q_0)^2$, $d_{2,1}=-2q_{21}^{-1}q_0\qnum{3}{q_0}$,
$\qnum{2}{p_3}=1+q_0$, and $\paar{\hatz _2}{t_2}=\paar{\hatz _2}{z_2}
(1-q_0^{-1})^{-1}q_{21}^{-2}q_0(1+q_0^2)^2(1+q_0)z_{2,2}$. Then
Corollary \ref{f-zi2=0} and Lemma \ref{l-t=0orinfty} give a
contradiction to (A0).

Suppose now that $a\ge 3$. Then
\begin{align*}
\frac{\paar{\hatz _{a-2}}{w_{a-1}}}{\paar{\hatz _{a-2}}{z_{a-2}}}=&
-(q_{11}^{-2}+1)\frac{q_{11}^4q_{21}^{-2}(1-q_{11}^{-5}q_{12}q_{21})
(1-q_{11}^{-2}{q_{12}q_{21}})}{1-q_{11}^{-1}}z_a.
\end{align*}
By (A8) and $(q_{12}q_{21})^{a+2}\not=1$ this implies that
$\paar{\hatz _{a-2}}{w_{a-1}}\not=0$ which is a contradiction to
(A7).

\mysubs
\label{ss-Serie6}
Finally we have to consider Nichols algebras where $q_{11}\in
R_{a+2}$, $q_{22}=-q_{11}^{-1}(q_{12}q_{21})^{-a-1}$,
and $(q_{12}q_{21})^{a+2}\not=1$ (cf.~Subsect.~\ref{ss-Serie4}).
One has
\begin{gather*}
\paar{y_2}{w_1}=\frac{(1+q^{-1}_{11}(q_{12}q_{21})^{-a})
(1-q_{11}(q_{12}q_{21})^{a+1})(1-(q_{12}q_{21})^{1-a})
}{q_{21}^2(1-(q_{12}q_{21})^{-a})}z_2.
\end{gather*}
If $a=1$ then by assumption (A3) one has $u_3=0$.
Equation $u_2=0$ is equivalent to $(q_{22}+1)(1-q_{12}q_{21}q_{22})=0$.
The cases $q_{22}=-1$ and $q_{12}q_{21}q_{22}=1$ were considered in
Subsection \ref{ss-Serie5} and \ref{ss-Serie4}, respectively.
Moreover, $u_3=0$, $u_2\not=0$ imply that
$\qnum{3}{q_{22}}(1-q_{12}q_{21}q_{22}^2)=0$. The algebra $\cB (V)$
with $q_{12}q_{21}q_{22}^2=1$ was already considered in Subsection
\ref{ss-Serie3} (exchange the generators $x_1$ and $x_2$). If
$\qnum{3}{q_{22}}=0$ then either $q_{11}=q_{22}^{-1}$ or
$q_{11}=q_{22}$. In the first case one has $q_0:=q_{11}q_{12}q_{21}$,
$q_0\in R_{12}$, $q_{11}=q_0^4$, $q_{22}=-q_0^2$.
Otherwise $\qnum{3}{-(q_{12}q_{21})^2}=0$ and
$q_{11}=q_{22}=-(q_{12}q_{21})^2$. These examples
appear in Theorem \ref{t-class}(5.1) and Theorem \ref{t-class}(5.2),
respectively.

If $a\ge 2$ then we must have $w_1=w_2=0$ by (A7).
For $(q_{12}q_{21})^a=q_{12}q_{21}$, $(q_{12}q_{21})^a=-q_{11}^{-1}$, and
$(q_{12}q_{21})^a=q_{11}^{-1}(q_{12}q_{21})^{-1}$ one computes
\begin{align*}
\frac{\paar{\hatz _1}{w_2}}{\paar{\hatz _1}{z_1}}=&
\frac{-q_{21}^{-2}(1+q_{11}^{-2})}{1-q_{11}^{-1}}
q_{11}^{-1}(q_{12}q_{21}-1)(q_{11}^3q_{12}q_{21}-1)z_3,\\
\frac{\paar{\hatz _1}{w_2}}{\paar{\hatz _1}{z_1}}=&
\frac{q_{21}^{-2}(q_{11}^{-3}(q_{12}q_{21})^{-1}-1)}{1+
q_{11}^{-4}(q_{12}q_{21})^{-1}}q_{11}^{-1}\qnum{2}{q_{11}^{-1}}
\qnum{2}{q_{11}^{-2}}\qnum{3}{-q_{11}^2q_{12}q_{21}}z_3,\\
\frac{\paar{\hatz _1}{w_2}}{\paar{\hatz _1}{z_1}}=&
\frac{-q_{21}^{-2}(1+q_{11}^{-3}(q_{12}q_{21})^{-2})}{
1+q_{11}^{-2}(q_{12}q_{21})^{-1}}
q_{11}^{-2}(1-q_{11}^3q_{12}q_{21})\qnum{3}{-q_{11}q_{12}q_{21}}z_3,
\end{align*}
respectively.
Let $a=2$. Then $w_1=0$, $q_{11}^2=-1$, and $(q_{12}q_{21})^4\not=1$ imply
that $q_{11}(q_{12}q_{21})^2=-1$ or $q_{11}=-(q_{12}q_{21})^3$.
In the first case one has $(q_{12}q_{21})^4=-1$, $q_{11}=(q_{12}q_{21})^2$,
$q_{22}=(q_{12}q_{21})^{-1}$ which was already considered in Subsection
\ref{ss-Serie5}. In the second case one has $\qnum{3}{-(q_{12}q_{21})^2}=0$,
$q_{11}=(q_{12}q_{21})^{-3}$, $q_{22}=-1$.
Then $\cB (V)$ appears in Theorem \ref{t-class}(4.6).

Assume that $a\ge 3$. Then (A8), $(q_{12}q_{21})^{a+2}\not=1$, and
$w_1=w_2=0$ imply that either Equations
$(q_{12}q_{21})^a=-q_{11}^{-1}$, $\qnum{3}{-q_{11}^2q_{12}q_{21}}=0$ or
$q_{11}(q_{12}q_{21})^{a+1}=1$, $q_{11}^3(q_{12}q_{21})^2=-1$ or
$q_{11}(q_{12}q_{21})^{a+1}=1$, $\qnum{3}{-q_{11}q_{12}q_{21}}=0$ hold.
Elementary computations show that in the first case there exist
precisely two
solutions satisfying the assumptions at the beginning of this subsection.
If $a=3$, $q_{12}q_{21}\in R_{30}$,
$q_{11}=-(q_{12}q_{21})^{-3}$, and $q_{22}=(q_{12}q_{21})^{-1}$ then
$\cB (V)$ appears in Theorem \ref{t-class}(2.7).
Otherwise $a=13$, $q_{12}q_{21}\in R_{30}$,
$q_{11}=(q_{12}q_{21})^2$, and $q_{22}=(q_{12}q_{21})^{-1}$. In this case
$\chi (z_9,z_9)=-(q_{12}q_{21})^5$ and
$d_{9,0}=2q_{21}^{-1}/(1-q_{12}^{-1}q_{21}^{-1})\not=0$
which is a contradiction to (A6).

In the second case, i.\,e.~if $q_{11}(q_{12}q_{21})^{a+1}=1$ and $q_{11}^3
(q_{12}q_{21})^2=-1$, (A8) and $(q_{12}q_{21})^{a+2}\not=1$ imply
that $a=3$, $q_{12}q_{21}\in R_{20}$,
$q_{11}=(q_{12}q_{21})^{-4}$, and $q_{22}=-1$.
Then $\cB (V)$ appears in Theorem \ref{t-class}(4.7).

Finally, Equations $q_{11}(q_{12}q_{21})^{a+1}=1$,
$\qnum{3}{-q_{11}q_{12}q_{21}}=0$ have to be considered.
One gets $q_{11}^3(q_{12}q_{21})^{3a+3}=1$ and
$(q_{12}q_{21})^{3(a+2)}=(-1)^a$. Hence $q_{11}^3(q_{12}q_{21})^{-3}=(-1)^a$
and together with $\qnum{3}{-q_{11}q_{12}q_{21}}=0$ one obtains
$q_{11}^6=(-1)^{a+1}$. By (A8) this gives $a=10$ and hence
$\qnum{3}{-(q_{12}q_{21})^{10}}=0$. Together with $1=q_{11}^{-a-2}=
(q_{12}q_{21})^{(a+1)(a+2)}=(q_{12}q_{21})^{132}$ this implies that
$q_{12}q_{21}$ is a $12^\mathrm{th}$ root of unity. The latter is a
contradiction to $(q_{12}q_{21})^{a+2}\not=1$.

Now all possible settings for the structure constants $q_{ij}$,
$i,j\in \{1,2\}$, are investigated. Thus the proof of Theorem \ref{t-class}
is finished.

\bibliography{quantum}
\bibliographystyle{mybib}

\end{document}